\documentclass[10pt,a4paper]{article}
\pdfoutput=1 

\usepackage{lineno,hyperref}
\modulolinenumbers[5]

\usepackage{geometry}
\usepackage{amssymb}
\usepackage{amsmath}
\usepackage{algorithm}     
\usepackage{algpseudocode} 
\usepackage{subfig}
\usepackage{csquotes}
\usepackage{tikz}
\usetikzlibrary{arrows}
\usepackage{mathtools}
\usepackage{bbm}
\usepackage{amsthm}
\usepackage{pgfplots}
\usepackage{float}
\usepgfplotslibrary{external}

\newcommand{\email}[1]{\hspace*{\stretch{1}}\emph{\texttt{#1}}}

\makeatletter
\def\blfootnote{\xdef\@thefnmark{$\star$}\@footnotetext}
\makeatother
\newenvironment{Authors}%
  {\begin{center}\begin{bfseries}}%
  {\end{bfseries}\end{center}}
\newenvironment{Addresses}%
  {\begin{flushleft}\begin{itshape}}%
  {\end{itshape}\end{flushleft}}

 \usepackage{fancyhdr}  
  
\fancypagestyle{plain}{
	\fancyhead{}
	\fancyhead[C]{\hfill submitted to Elsevier, May 2021}
}
  
\newtheorem{theorem2}{Theorem}[section]

\newtheorem{remark}[theorem2]{Remark}

  \newcommand{\vertiii}[1]{{\left\vert\kern-0.25ex\left\vert\kern-0.25ex\left\vert #1 
    \right\vert\kern-0.25ex\right\vert\kern-0.25ex\right\vert}}

\begin{document}

\thispagestyle{plain}

\title{Registration-based model reduction of parameterized two-dimensional conservation laws.}
 \date{}
 
 \maketitle
\vspace{-50pt} 
 
\begin{Authors}
Andrea Ferrero$^{1}$,
Tommaso Taddei$^{2}$,
Lei Zhang$^{2}$
\end{Authors}

\begin{Addresses}
$^1$
Department of Mechanical and Aerospace Engineering, Politecnico di Torino, Corso Duca degli Abruzzi 24, 10129 Torino, Italy 
 \email{andrea{\_}ferrero@polito.it} \\[3mm]
$^2$
IMB, UMR 5251, Univ. Bordeaux;  33400, Talence, France.
Inria Bordeaux Sud-Ouest, Team MEMPHIS;  33400, Talence, France, \email{tommaso.taddei@inria.fr,lei.a.zhang@inria.fr} \\
\end{Addresses}

\begin{abstract}
We propose a nonlinear registration-based model reduction procedure for rapid and reliable solution of  parameterized two-dimensional steady conservation laws.
This class of problems is challenging for model reduction techniques due to the presence of nonlinear terms in the equations and also due to the presence of parameter-dependent discontinuities that cannot be adequately represented through linear approximation spaces.
Our approach builds on 
a general (i.e., independent of the underlying equation) registration procedure for the computation of a mapping $\Phi$ that tracks moving features of the solution field and on an hyper-reduced least-squares Petrov-Galerkin  reduced-order model  for the rapid and reliable computation of the solution coefficients.
The contributions of this work are twofold.
First, we investigate the application of registration-based methods to two-dimensional hyperbolic systems.
Second, we propose a multi-fidelity approach to reduce the offline costs associated with the construction of the parameterized mapping and the reduced-order model.
We discuss the application to an inviscid supersonic flow past a parameterized bump, to illustrate the many features of our method and to demonstrate its effectiveness.
\end{abstract}

\medskip

 
\section{Introduction}
\label{sec:intro}

\subsection{Model order reduction for steady conservation laws}
\label{sec:intro_a}
Despite the recent advances in high-performance computing and numerical analysis,  approximation of the solution to fluid problems remains a formidable task that requires extensive computational resources. The lack of fast and reliable computational fluid dynamics (CFD) solvers limits the use of high-fidelity (hf) simulations to perform extensive parametric studies in science and engineering. Parameterized model order reduction (pMOR) aims at constructing a low-dimensional  surrogate (or reduced-order) model (ROM) over a range of parameters, and ultimately speed up parametric studies. The goal of this paper is to develop a nonlinear registration-based MOR procedure for steady two-dimensional conservation laws and to demonstrate its effectiveness for applications in aerodynamics.

We denote by ${\mu}$ the vector of model parameters in the parameter region $\mathcal{P} \subset \mathbb{R}^P$; we denote by $\Omega \subset \mathbb{R}^2$ the computational domain --- to simplify the presentation, in the introduction we assume that the domain does not depend on the parameters; however, in the numerical examples, we shall consider the case of parameterized geometries.  We denote by ${U}: \Omega \times \mathcal{P} \to \mathbb{R}^D$ the parametric solution field satisfying the conservation law:
\begin{equation}
\label{eq:conservation_law}
\nabla \cdot {F}_{\mu}( {U}_{\mu} ) = {S}_{\mu}({U}_{\mu}) \quad
{\rm in} \; \Omega,
\end{equation}
where ${F}: \mathbb{R}^D \times \mathcal{P} \to \mathbb{R}^{D,2}$ is the physical flux and 
${S}: \mathbb{R}^D \times \mathcal{P} \to \mathbb{R}^{D}$ is the source term. The problem is completed with suitable boundary conditions that depend on the number of incoming characteristics. We denote by $\mathcal{M} : = \{ {U}_{\mu} : \mu \in \mathcal{P}   \}$ the solution manifold associated with \eqref{eq:conservation_law}. We further define the Hilbert space $\mathcal{X} = [L^2(\Omega)]^D$, endowed with the inner product $(\cdot, \cdot)$ and the induced norm
$\| \cdot \|:=\sqrt{(\cdot,\cdot)}$, such that
$({w}, {v}) = \int_{\Omega} {w} \cdot  {v} \,  d{x}$ for all ${w}, {v} \in \mathcal{X}$.

We introduce the finite element (FE) mesh 
 $\mathcal{T}_{\rm hf}: = \left(  \{ {x}_{j}^{\rm hf} \}_{j=1}^{N_{\rm hf,v}}, \texttt{T} \right)$ where 
 $ \{ {x}_{j}^{\rm hf} \}_j \subset \overline{\Omega}$ are the nodes of the mesh and $\texttt{T} \in \mathbb{N}^{n_{\rm lp}, N_{\rm e}}$ is the connectivity matrix, where  $n_{\rm lp}$ is the number of degrees of freedom in each element and $N_{\rm e}$ is the total number of elements.  
 We denote by $\mathcal{X}_{\rm hf} \subset \mathcal{X}$ a FE discretization associated with $\mathcal{T}_{\rm hf}$ and we set $N_{\rm hf} = {\rm dim} ( \mathcal{X}_{\rm hf}  )$. Given ${w} \in \mathcal{X}_{\rm hf}$, we denote by ${\mathbf{w}} \in \mathbb{R}^{N_{\rm hf}}$ the vector representation of ${w}$ with respect to a suitable basis: note that the pair mesh-coefficients $(\mathcal{T}_{\rm hf}, \;  {\mathbf{w}}  )$ uniquely identifies the field ${w} \in \mathcal{X}_{\rm hf}$. Finally, we denote by ${U}_{\mu}^{\rm hf}\in \mathcal{X}_{\rm hf}$ the hf estimate of the solution ${U}_{\mu} \in \mathcal{X}$ to \eqref{eq:conservation_law} for a given ${\mu} \in \mathcal{P}$.
 
Hyperbolic problems with moving fronts are extremely challenging for state-of-the-art model reduction procedures.
First, the vast majority of MOR methods rely on 
linear approximations: 
as shown in several studies (e.g.,
\cite{ohlberger2015reduced}), linear methods  are fundamentally ill-suited to deal with parameter-dependent sharp gradients that naturally  arise in the  solutions to hyperbolic conservation laws.
 Another major issue concerns the construction of accurate meshes for parametric studies. For advection-dominated problems, adaptive mesh refinement (AMR) is of paramount importance to reduce the size of the mesh required to achieve a given accuracy. However, if parametric variations strongly affect the location of sharp-gradient regions, AMR should be applied to each system configuration and will lead to hf discretizations of intractable size. Effective MOR procedures for conservation laws should thus embed an effective parametric AMR strategy to track moving structures. 
 
\subsection{Registration methods for parameterized problems}
\label{sec:intro_b}

Registration-based  (or Lagrangian) methods 
for pMOR
(e.g.,  \cite{iollo2014advection,ohlberger2013nonlinear,sarna2020data,taddei2015reduced,taddei2020registration}
) rely on the introduction of a parametric mapping ${\Phi}: \Omega \times \mathcal{P} \to \Omega$ such that (i) ${\Phi}_{\mu}$ is a bijection from $\Omega$ in itself for all ${\mu} \in \mathcal{P}$, and 
(ii) the mapped manifold $\widetilde{\mathcal{M}} = \{ 
{U}_{\mu} \circ {\Phi}_{\mu}:  {\mu} \in \mathcal{P}
\}$ is more amenable for linear compression methods. In the FE framework, or equivalently in the finite volume context, this corresponds to considering approximations of the form
\begin{equation}
\label{eq:registration_abstract}
{\mu} \in \mathcal{P} \mapsto
\left(
{\Phi}_{\mu}( \mathcal{T}_{\rm hf}   ), \;
\widehat{{\mathbf{U}}}_{\mu} = \mathbf{Z} \;  \widehat{\boldsymbol{\alpha}}_{\mu}
\right),
\quad
{\rm with} \;
{\Phi}_{\mu}( \mathcal{T}_{\rm hf}   ): = \left(  \{ {\Phi}_{\mu}({x}_{j}^{\rm hf} )\}_{j=1}^{N_{\rm hf,v}}, \texttt{T} \right),
\;\;\;
\mathbf{Z} \in \mathbb{R}^{N_{\rm hf}, N}.
\end{equation}
Note that the mapped mesh 
${\Phi}_{\mu}( \mathcal{T}_{\rm hf}   )$ shares with 
$\mathcal{T}_{\rm hf} $ the same connectivity matrix, while 
$\widehat{{\mathbf{U}}}_{\mu} = \mathbf{Z} \;  \widehat{\boldsymbol{\alpha}}_{\mu}$ can be viewed as an approximation of ${U}_{\mu}$ if paired with the mesh 
${\Phi}_{\mu}( \mathcal{T}_{\rm hf}   ) $, or as an approximation of 
 ${U}_{\mu}\circ {\Phi}_{\mu}$ if paired with the mesh $\mathcal{T}_{\rm hf}$.
 
Several features of registration methods are attractive for applications to hyperbolic problems with moving fronts.
First, registration methods are effective to track sharp gradients of the solution field, and ultimately improve performance of linear compression methods in the reference configuration and also reduce the 
size of the hf mesh required for a given accuracy.
Second,  after having built the mapping 
 ${\Phi}$, 
Lagrangian methods reduce to linear methods in parameterized domains: this class of methods has been widely studied in the MOR literature 
(see the reviews
\cite{lassila2014model,rozza2021basic} and also
\cite{taddei2020discretize}) and is now well-understood.
  In particular, we can rely on standard 
 training algorithms to build $\mathbf{Z}$ in \eqref{eq:registration_abstract}
 -- in particular,   proper orthogonal decomposition (POD, \cite{berkooz1993proper,volkwein2011model}) and the weak-Greedy algorithm \cite{rozza2007reduced} --- and on effective hyper-reduced projection-based  techniques to compute the solution  coefficients $\widehat{{\boldsymbol{\alpha}}}_{\mu}$.

In this work, we consider  the registration procedure first proposed in \cite{taddei2020registration} and then extended in  \cite{taddei2021space,taddei2021registration} to generate the mapping; then, similarly to 
\cite{taddei2021space}, we rely on a projection-based least-squares Petrov-Galerkin (LSPG, \cite{carlberg2011efficient,carlberg2017galerkin}) formulation with elementwise empirical quadrature (EQ, \cite{farhat2015structure,yano2019discontinuous})  to estimate the coefficients 
$\widehat{\boldsymbol{\alpha}}_{\mu}$ for any new value of the parameters.
Furthermore, we rely on the \emph{discretize-then-map} framework 
(cf. 
\cite{dal2019hyper,taddei2020discretize,washabaugh2016use})
to deal with geometry variations. 
 The contribution of the paper is twofold.
\begin{itemize}
\item
We show performance of registration-based model reduction for a representative problem in aerodynamics  with shocks: we discuss performance of registration, and we also address the combination with projection-based MOR techniques. In particular, we investigate  in detail the offline-online computational decomposition and we also comment on hyper-reduction, which is key for online efficiency.
\item
We present work toward the implementation of a multi-fidelity approach for registration-based model reduction. As explained in 
\cite{taddei2020registration,taddei2021space,taddei2021registration},
a major issue of our registration procedure is the need for extensive explorations of the parameter domain:  in this work, we show that we can rely on a significantly less  accurate hf discretization to generate the snapshots used for registration and ultimately greatly reduce the cost of offline training. In the numerical results, we further show that multi-fidelity training might help reduce the size of the hf discretization required to properly track moving features --- in effect,    spatio-parameter mesh adaptivity.
\end{itemize}

The outline of the paper is as follows.
In section \ref{sec:model_problem}, we introduce the model problem; in section \ref{sec:methods}, we present the methodology: first, we introduce the registration algorithm proposed in \cite{taddei2021registration}, then, we discuss the  projection-based scheme and finally we present the offline/online computational decomposition based on a two-fidelity sampling.
In section \ref{sec:numerics}, we present extensive numerical investigations to illustrate the performance of our proposal.
In the remainder of this section, we discuss relation to previous works (cf. section \ref{sec:prior_work}),  we briefly comment on the many  nonlinear approximation methods 
appeared in the literature to better clarify the interest for  registration-based methods (cf. section \ref{sec:review_nonlinear}), and we present relevant 
notation (cf. section  \ref{sec:notation}).

\subsection{Relation to previous works}
\label{sec:prior_work}

Several authors have applied MOR techniques to aerodynamics problems including inviscid flows: we refer to \cite{yano2021model} for a review; we further refer to the   early work by Zimmermann et al, \cite{franz2014interpolation} and to  the more recent work  by Carlberg et al, 
\cite{blonigan2021model} for application to aerodynamics  of techniques based on nonlinear  approximations.
Simultaneous adaptivity in space --- via AMR --- and in parameter --- via Greedy sampling --- has been considered by Yano in \cite{yano2018reduced} and more recently in \cite{sleeman2021goal}.
Methods in \cite{sleeman2021goal,yano2018reduced} rely on $h$-refinement to adapt the spatial mesh,  while we  exploit a solution-aware parameterized  
mapping to deform the mesh without changing its topology ($r$-adaptivity): we thus envision that the two strategies  might be combined with mutual benefits.

Multifidelity methods have been extensively studied in the MOR literature: we refer to 
\cite{peherstorfer2018survey} for a thorough review and also  to the more recent work by Kast et al. \cite{kast2020non}. As explicitly stated in section \ref{sec:intro_b}, the present study offers a proof of concept of the application of multifidelity schemes in combination with registration methods; it also shows the importance of multifidelity schemes for spatio-parameter adaptivity.  

As discussed in 
\cite{taddei2020registration,taddei2021space,taddei2021registration},
the fundamental building block of our registration procedure is a nonlinear non-convex optimization 
statement for the computation of the mapping $\Phi$ for the parameters in the training set. 
Our optimization statement minimizes an $L^2$ reconstruction error 
 plus a number of terms that control the smoothness of the map and the mesh distortion: minimization of the $L^2$ reconstruction error has been previously considered in  several works (e.g., 
\cite{mojgani2020physics,reiss2020optimization,rim2018transport,sarna2020model}); on the other hand, penalization of 
mesh distortion has been considered in 
\cite{zahr2018optimization} in a related context.

For completeness, as already discussed  in 
\cite{taddei2021registration}, we remark that
registration-based methods  are tightly linked to a number of techniques in  related fields.
First, registration is central in image processing: in this field, registration refers to the process of transforming different sets of data into one coordinate system, \cite{zitova2003image}.
In computational mechanics,
Persson and Zahr have proposed in \cite{zahr2018optimization}  an $r$-adaptive optimization-based 
high-order discretization method 
 to deal with shocks/sharp gradients of the solutions to advection-dominated problems.
In uncertainty quantification, several authors (see, e.g., \cite{marzouk2016sampling}) have proposed measure transport approaches to sampling: transport maps are used to ``push forward'' samples from a reference configuration and ultimately facilitate sampling from non-Gaussian distributions.
Finally, the notion of registration is also at the core of diffeomorphic dimensionality reduction (\cite{walder2008diffeomorphic}) in the field of machine learning.

\subsection{Methods based on nonlinear approximations: expressivity and learnability}
\label{sec:review_nonlinear}

In recent years, there has been   a growing interest in nonlinear model reduction techniques, particularly for CFD applications.
A first class of methods relies on adaptive partitioning of the parameter domain, \cite{eftang2010hp}.
Another class of methods relies on online basis update and/or refinement:
relevant works that fit in this category might rely on  low-rank updates (e.g.,  \cite{carlberg2015adaptive,etter2020online,peherstorfer2020model}), 
or might rely on 
Grasmannian learning to construct 
parameter-dependent reduced-order bases  
\cite{amsallem2008interpolation,zimmermann2018geometric}.
A third  class of methods relies on the introduction in the offline/online workflow of a preprocessing stage to   reformulate the problem  in a form that is more amenable for linear approximations: representative methods in this category are the 
the approach in  \cite{gerbeau2014approximated}
based on approximate Lax pairs, and the method of freezing in \cite{ohlberger2013nonlinear}.
We remark that such preprocessing stage might be performed once during the offline stage, at the beginning of the online stage for any new $\mu \in \mathcal{P}$, or at each time step in combination with a suitable time-marching scheme.
A fourth class of methods considers directly  nonlinear approximations in combination with specialized methods to compute the solution during the online stage: to provide concrete references, we refer to the approaches based on convolutional autoencoders, \cite{fresca2021comprehensive,kashima2016nonlinear,kim2020efficient,lee2020model},
and to the  approach in  \cite{ehrlacher2020nonlinear} 
based on optimal transport and nonlinear interpolation.
As explained below (cf. \eqref{eq:abstract_class_lagrangian}), Lagrangian methods lead to predictions $\widehat{U}$ that are linear in the solution coefficients $\widehat{\boldsymbol{\alpha}}_{\mu}$ and nonlinear in the mapping coefficients $\widehat{\mathbf{a}}_{\mu}$: 
depending on the way mapping coefficients are computed, Lagrangian methods fit in the third category (e.g., \cite{taddei2020registration,taddei2021space} and this work) or in the fourth category (e.g., 
\cite{mojgani2017arbitrary}).

To analyze the many nonlinear proposals and ultimately perform an informed decision for the specific problem of interest, we shall interpret pMOR techniques as the combination of two  fundamental blocks: a low-rank parameter-independent operator 
$\texttt{Z}: \mathcal{A}\subset  \mathbb{R}^Q  \to \mathcal{X}$ and a ROM for the reduced coefficients  $\widehat{\boldsymbol{\beta}}: \mathcal{P} \to \mathcal{A}$.
To build $\texttt{Z}$, we first identify a class of approximations (see \eqref{eq:abstract_class} below) and then we proceed to use offline data to identify the proper (quasi-optimal) approximation within that class;  after having built $\texttt{Z}$, we rely on intrusive (projection-based) or non-intrusive (data-fitted) methods to rapidly find the coefficients $\widehat{\boldsymbol{\beta}}_{\mu} \in \mathcal{A}$ for any new value of the parameters in $\mathcal{P}$.
Examples of approximation classes include the aforementioned linear methods, Lagrangian methods, convolutional methods, and transported methods. 
\begin{subequations}
\label{eq:abstract_class}
\begin{itemize}
\item
Linear methods can be written as
\begin{equation}
\label{eq:abstract_class_linear}
\widehat{U}_{\mu} = 
\texttt{Z}(  \widehat{\boldsymbol{\beta}} _{\mu}= \widehat{\boldsymbol{\alpha}}_{\mu}  )
=
\sum_{n=1}^N  (\widehat{\boldsymbol{\alpha}}_{\mu}   )_n
\zeta_n,
\end{equation}
with $N=Q$,
$\mathcal{A} = \mathbb{R}^N$, 
  and $\zeta_1,\ldots,\zeta_N \in \mathcal{X}$. 
\item
Lagrangian (or registration-based) methods can be written as
\begin{equation}
\label{eq:abstract_class_lagrangian}
\widehat{U}_{\mu} = 
\texttt{Z}(  \widehat{\boldsymbol{\beta}}_{\mu} = [ \widehat{\boldsymbol{\alpha}}_{\mu}, \widehat{\mathbf{a}}_{\mu}]  )
=
\sum_{n=1}^N  (\widehat{\boldsymbol{\alpha}}_{\mu} )_n
\zeta_n \circ  \texttt{N}(  \widehat{\mathbf{a}}_{\mu} )^{-1}
\end{equation}
where $\zeta_1,\ldots,\zeta_N \in \mathcal{X}$, 
$\mathcal{A} = \mathbb{R}^N \times \mathcal{A}_{\rm bj}$, 
$\texttt{N}:  \mathbb{R}^M \to {\rm Lip}(\Omega; \mathbb{R}^2)$ such that
$\texttt{N}(  \mathbf{a}  )$ is a bijection in $\Omega$ for all $\mathbf{a} \in \mathcal{A}_{\rm bj}$,
$Q=  N+M$. 
\item
Convolutional approximations (\cite{fresca2021comprehensive,kashima2016nonlinear,kim2020efficient,lee2020model})
 with $L>0$ layers can be stated as
\begin{equation}
\label{eq:abstract_class_convolutional}
\texttt{Z}(  \widehat{\boldsymbol{\beta}}_{\mu} = [\widehat{\boldsymbol{\alpha}}_{1,\mu}, \ldots,  \widehat{\boldsymbol{\alpha}}_{L,\mu}])
=
\texttt{N}_L\left(
\texttt{N}_{L-1}\left(\ldots, \;
\boldsymbol{\alpha}_{L-1,\mu}
\right),
\boldsymbol{\alpha}_{L,\mu}
\right)  
\end{equation}
where
$\texttt{N}_{\ell}: \mathbb{R}^{D_{\ell}} \times\mathbb{R}^{N_{\ell}} \to \mathbb{R}^{D_{\ell+1}}$
with $D_{1}=2$ (number of spatial dimensions) and 
$D_{L+1}=D$  (number of state variables),
$Q=\sum_{\ell=1}^L N_{\ell}$ and
$\mathcal{A} = \mathbb{R}^Q$.
\item
Finally, 
transported (or transformed) snapshot methods
(\cite{cagniart2019model,nair2018transported,reiss2018shifted,welper2017interpolation})
with $N>0$ terms can be stated as
\begin{equation}
\label{eq:abstract_class_transported}
\texttt{Z}(  \widehat{\boldsymbol{\beta}}_{\mu} = [ \widehat{\boldsymbol{\alpha}}_{\mu}, \widehat{\mathbf{a}}_{1,\mu},\ldots,\widehat{\mathbf{a}}_{N,\mu}]  )
=
\sum_{n=1}^N  (\widehat{\boldsymbol{\alpha}}_{\mu} )_n
\zeta_n \circ \texttt{N}_n(  \widehat{\mathbf{a}}_{n,\mu}  )
\end{equation}
where 
$\texttt{N}_1, \ldots, \texttt{N}_N: \mathbb{R}^M \to {\rm Lip} (\Omega: \mathbb{R}^2)$,
$\zeta_1,\ldots,\zeta_N \in \mathcal{X}_{\rm ext} : = \{ v\in L^2(\mathbb{R}^2): v|_{\Omega} \in \mathcal{X} \}$.
\end{itemize}
Note that, while in Lagrangian methods we require that $\texttt{N}$ is bijective, transformed methods  do not explicitly require bijectivity of 
$\texttt{N}_1, \ldots,\texttt{N}_N$.   Note also that linear methods are a subset of Lagrangian methods --- in the sense that they reduce to linear methods for  $\texttt{N} = \texttt{id}$, 
$\texttt{id}(x) = x$. Similarly, Lagrangian methods are a subset of convolutional and transported methods.
\end{subequations}

The  choice of the class of approximations should be a compromise between \emph{expressivity}  and \emph{learnability}.
In statistical learning, expressivity (or expressive power) of a network refers to the approximation properties for a given class of functions, \cite{guhring2020expressivity}.
Given the class of approximations $\mathcal{C} \subset C(\mathcal{A}; \mathcal{X} )$ for some $Q> 0$ ---
 $C(\mathcal{A}; \mathcal{X} )$ is the space of continuous applications from $\mathcal{A} \subset \mathbb{R}^Q$ to  $\mathcal{X}$
 --- we measure the expressivity of $\mathcal{C}$ for $\mathcal{M}$ in terms of the nonlinear width (\cite{dung1996nonlinear}):
 \begin{equation}
 \label{eq:nonlinear_width}
\inf_{ \texttt{Z} \in \mathcal{C}      } 
 \sup_{{w} \in \mathcal{M}} \;
 \inf_{\boldsymbol{\beta} \in \mathcal{A} }
 \;
\| \texttt{Z} (\boldsymbol{\beta}) - {w}  \|.
 \end{equation}
 On the other hand, learnability depends on two distinct factors: (i) the performance of available training algorithms to identify an approximation map $\texttt{Z}$ in $\mathcal{C}$
that approximately realizes the optimum of \eqref{eq:nonlinear_width}; and
 (ii) the performance  of available  methods to rapidly and reliably compute the coefficients $\widehat{\boldsymbol{\beta}}_{\mu}$ during the online stage.
Note that the training algorithm in (i) 
 is fed with a finite set of snapshots from $\mathcal{M}$: due to the large cost of hf CFD simulations, reduction of the number of required offline  simulations is key for  practical applications.
 
 Since expressivity depends on the particular manifold of interest, while learnability depends on the PDE model under consideration, it seems difficult to offer a definitive answer concerning the optimal choice of the approximation class $\mathcal{C}$. The aim of this work is to show that Lagrangian approximations have high expressive power for a representative problem in aerodynamics and that they can be learned effectively based on sparse datasets: further theoretical and numerical investigations are needed to clarify the scope of the present   class of methods and ultimately offer guidelines for the choice of the class of approximations.

\subsection{Notation}
\label{sec:notation}

 We estimate the solution to \eqref{eq:conservation_law} using a nodal-based discontinuous Galerkin (DG) finite element (FE) discretization of degree \texttt{p}.
Similarly to \cite{taddei2021registration}, we resort to
 a FE isoparametric discretization.
We define the reference element 
$\widehat{\texttt{D}} = \{ {X}\in [0,1]^2: \sum_{d=1}^2 X_d<  1  \}$  and 
  the Lagrangian basis $\{ \ell_i \}_{i=1}^{n_{\rm lp}}$  of the polynomial space $\mathbb{P}_{\texttt{p}}(\widehat{\texttt{D}})$ associated with the nodes 
 $\{  {X}_i  \}_{i=1}^{n_{\rm lp}}$;
 then, 
recalling the definition of $\mathcal{T}_{\rm hf}$ in section \ref{sec:intro_a}, 
  we define the  elemental mappings $\{ {\Psi}_k^{\rm hf} \}_{k=1}^{N_{\rm e}}$ such that
 \begin{equation}
\label{eq:psi_mapping}
 {\Psi}_k^{\rm hf}( {{X}})
 =
 \sum_{i=1}^{n_{\rm lp}} \;
x_{\texttt{T}_{i,k}}^{\rm hf}  \; \ell_i({X}),
\end{equation}
and  the elements of the mesh
$\{ \texttt{D}_k : =  {\Psi}_k( \widehat{\texttt{D}}   )\}_k$.
We further define the basis functions
$\ell_{i,k}: \Omega \to \mathbb{R}$ such that
$\ell_{i,k}(x) = 0$ for all $x\notin \texttt{D}_k$ and 
$\ell_{i,k} =  \ell_i \circ  {\Psi}_k^{-1}(x)$
for $x\in \texttt{D}_k$, $i=1,\ldots,n_{\rm lp}$, $k=1,\ldots,N_{\rm e}$.

We define the FE space $\mathcal{X}_{\rm hf} = {\rm span}\{
\ell_{i,k} {e}_d:
i=1,\ldots,n_{\rm lp}, k=1,\ldots,N_{\rm e},
d=1,\ldots,D \}$ where $e_1,\ldots,e_D$ are the canonical basis of $\mathbb{R}^D$. Given ${w} \in \mathcal{X}_{\rm hf}$, we denote by ${\mathbf{w}} \in \mathbb{R}^{N_{\rm hf}}$, $N_{\rm hf} = n_{\rm lp} \cdot N_{\rm e} \cdot D$, the corresponding vector of coefficients such that
\begin{equation}
\label{eq:vector2field}
{w}({x})
\,  = \, 
  \sum_{k=1}^{N_{\rm e}} \, 
\sum_{i=1}^{n_{\rm lp}} \, 
\sum_{d=1}^{D} \,  
\left( \mathbf{w} \right)_{i + n_{\rm lp} (k-1) + n_{\rm lp} N_{\rm e} (d-1)} \; \ell_{i,k}({x}) \; {e}_d,
\qquad \forall \; {x} \in \Omega.
\end{equation}
 Note that
\eqref{eq:vector2field} introduces an isomorphism 
between $\mathbb{R}^{N_{\rm hf}}$ and $\mathcal{X}_{\rm hf}$.
 Following the discussion in  \cite{taddei2020discretize}, we can extend the previous definitions to the mapped mesh and mapped FE space. We omit the details.
 
In view of the FE approximation, it is important that the deformed mesh 
$\Phi_{\mu}(\mathcal{T}_{\rm hf})$ (cf. \eqref{eq:registration_abstract}) does not have inverted elements. In this respect, we say that the mapping 
$\Phi: \Omega \times \mathcal{P} \to \mathbb{R}^2$ is bijective with respect to $\mathcal{T}_{\rm hf}$
(\emph{discrete bijectivity}, \cite[Definition 2.2]{taddei2021registration})
 if the elemental mappings of the deformed mesh are invertible.

\section{Model problem}
\label{sec:model_problem}

We consider the problem of approximating the solution  to  the parameterized  compressible  Euler equations.
The compressible Euler equations are a widely-used model to study aerodynamic flows: we refer to 
\cite{toro2013riemann} for a thorough discussion; we here consider the non-dimensional form of the equations.
We denote by $\rho$ the density of the fluid, by ${u}=[u_1,u_2]$ the velocity field,  by $E$ the total energy and by $p$ the pressure; we further define the vector of conserved variables ${U} = [ \rho, \rho {u}, E]: \Omega \to \mathbb{R}^{D=4}$. In this work, we consider the case of  ideal gases for which we have the following relationship between pressure and conserved variables ${U}$:
\begin{subequations}
\label{eq:euler_equations}
\begin{equation}
p = (\gamma - 1) \left( E - \frac{1}{2} \rho \|  {u} \|_2^2 \right),
\end{equation}
where $\gamma$ is  the ratio  of  specific  heats, which is here set equal to $\gamma=1.4$.
We further introduce the speed of sound $a$ and the Mach number ${\rm Ma}$ with respect to the channel axis such that 
\begin{equation}
\label{eq:a_mach}
a = \sqrt{\gamma \frac{p}{\rho}},
\quad
{\rm Ma} = \frac{u_1}{a}.
\end{equation}
Finally, we introduce the Euler physical flux and source term:
\begin{equation}
{F} ({U})
=\left[
\begin{array}{l}
\rho {u}^T \\
\rho {u} \, {u}^T
+ p \mathbbm{1}  \\
{u}^T (E+p) \\
\end{array}
\right],
\quad
{S} ({U}) = {0}.
\end{equation}
\end{subequations}

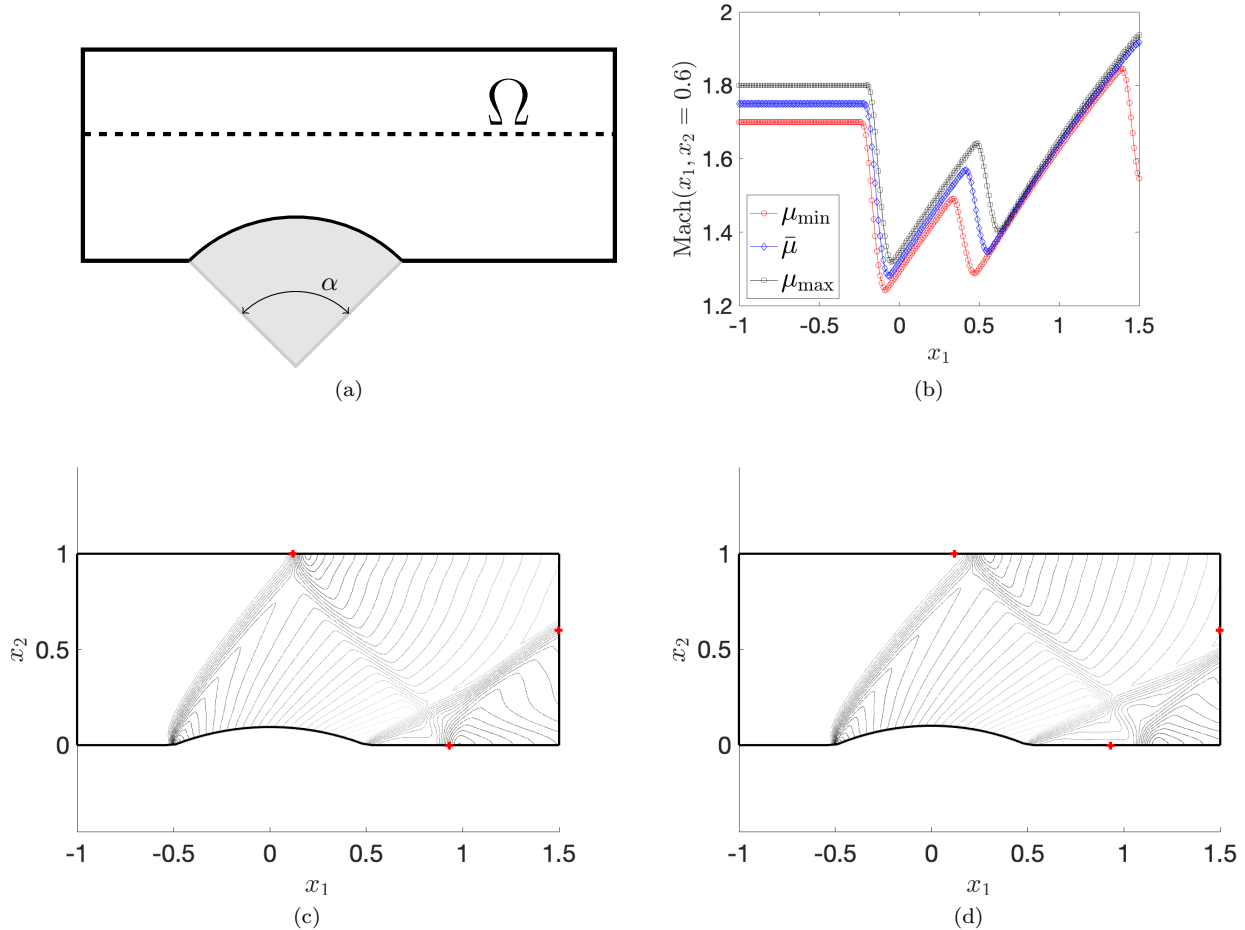
\begin{figure}[h!]
\centering

\subfloat[ ]{
\begin{tikzpicture}[scale=2.8]
\linethickness{0.3 mm}
\linethickness{0.3 mm}

\draw[ultra thick]  (-0.5,0)--(-1,0)--(-1,1)--(1.5,1)--(1.5,0)--(0.5,0);

\draw[ultra thick,dashed]  (-1,0.6)--(1.5,0.6);

\draw[very thick] (0.5,0)  arc(45:135:0.7071);


\draw[very thick, fill=gray, opacity=0.2] (0,-0.5) --(0.5,0) arc(45:135: 0.7071) -- cycle;

\draw [<->] (0.25,-0.25)  arc (45:135:0.3536) node[above,pos=.2]{$\alpha$};

\coordinate [label={above:  {\Huge {$\Omega$}}}] (E) at (1, 0.6) ;
\end{tikzpicture}
}
~~~
\subfloat[ ]{
\includegraphics[width=0.4\textwidth]
{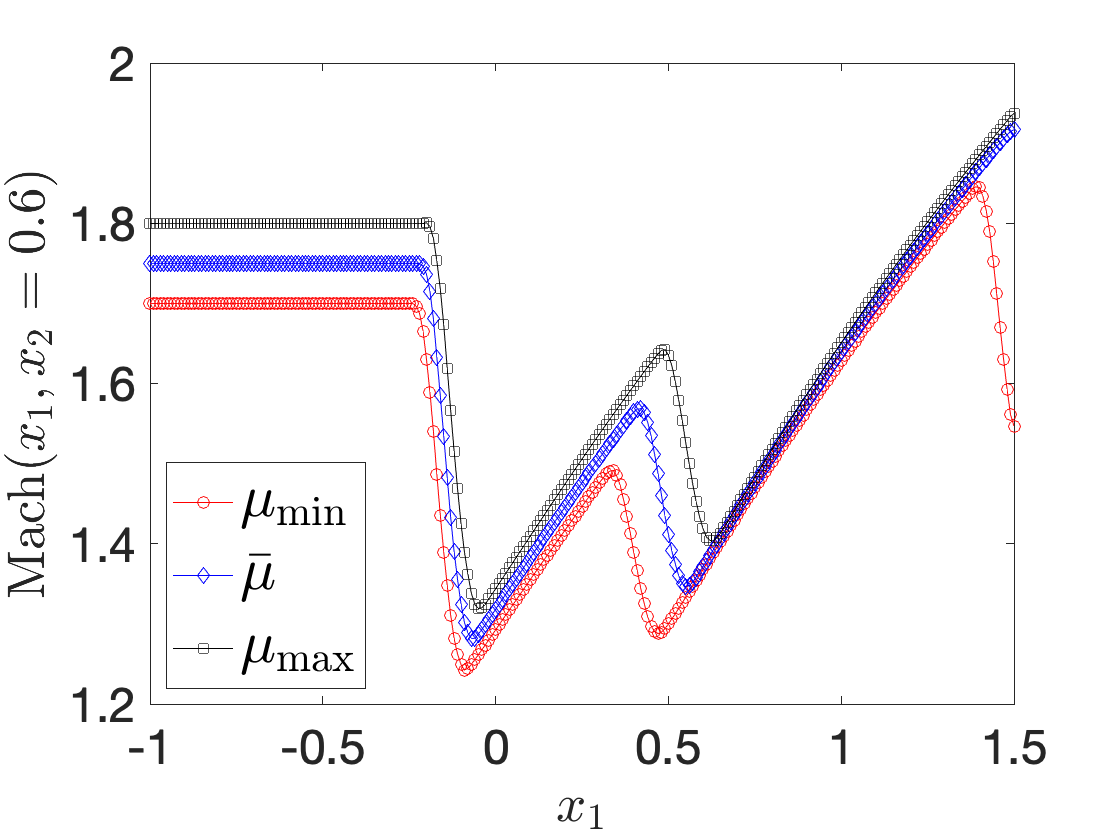}
}

\subfloat[ ]{
\includegraphics[width=0.48\textwidth]
{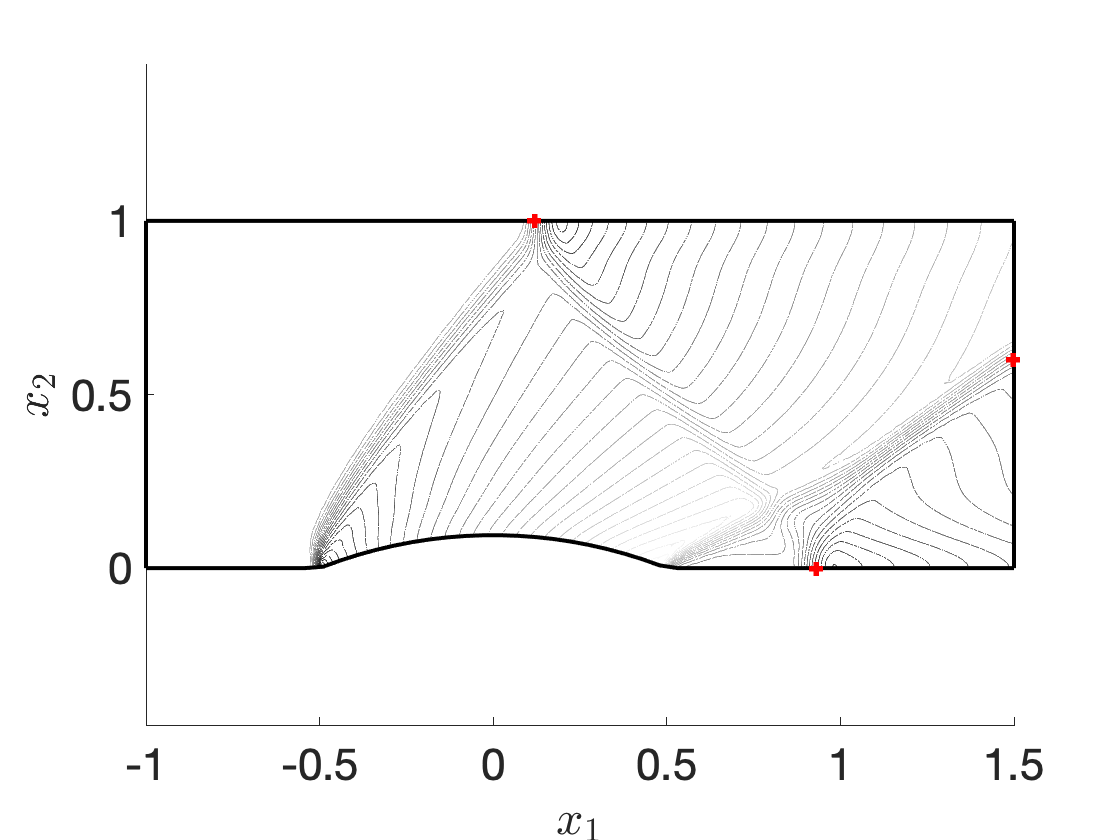}
}
~~~
\subfloat[ ]{
\includegraphics[width=0.48\textwidth]
{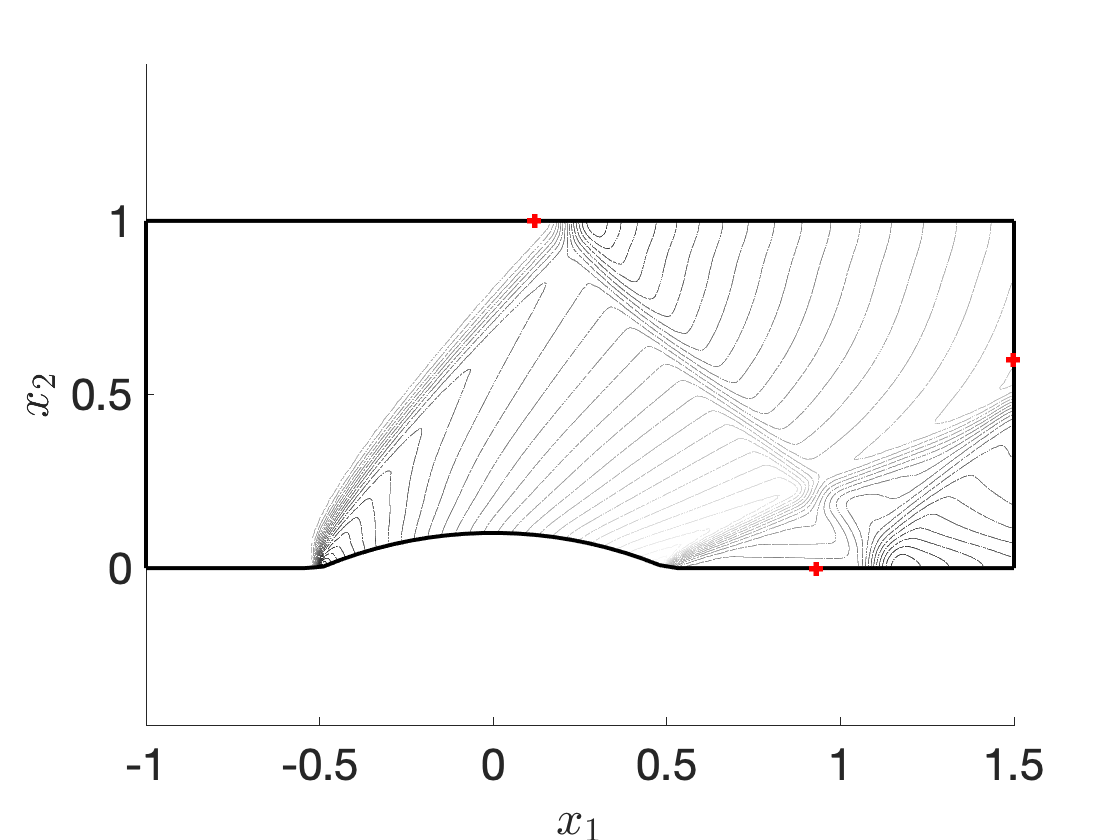}
}

\caption{flow past a circular bump. (a) geometric configuration. (b) horizontal  slices of the Mach number at $x_2=0.6$ for ${\mu}_{\rm min} = [0.75,1.7]$, 
 $\bar{{\mu}} = [0.775,1.75]$
 ${\mu}_{\rm max} = [0.8,1.8]$.
 (c)-(d) contour lines of the Mach number for 
 ${\mu}_{\rm min}$ and 
 ${\mu}_{\rm max}$.
}
\label{fig:bumpflow_vis}
\end{figure}

We consider a parametric channel flow past a circular bump:
the parameters are the free-stream Mach number ${\rm Ma}_{\infty}$ and the central angle $\alpha$ associated with the bump --- cf.  Figures \ref{fig:bumpflow_vis}(a), 
\begin{equation}
\label{eq:parameters}
{\mu} = [\alpha, {\rm Ma}_{\infty}  ]
\in \mathcal{P} = [0.75,0.8]\times [1.7,1.8].
\end{equation}
The horizontal length of the bump and the height of the channel are set to one. We impose wall conditions at the lower and upper boundaries, transmissive boundary conditions at the outflow and we set 
${U} = {U}_{\infty}$ at the inflow with 
$$
\rho_{\infty} = \frac{p_{\infty}}{T_{\infty}}, \quad
{u}_{\infty} =
 \sqrt{\gamma T_{\infty}} 
\left[
\begin{array}{l}
{\rm Ma}_{\infty} \\0\\
\end{array}
\right],
\quad
p_{\infty} = \frac{1}{ ( 1 +\frac{\gamma-1}{2} {\rm Ma}_{\infty}^2 )^{\frac{\gamma}{\gamma-1}}},
\quad
{\rm and}
\quad
T_{\infty} = \frac{1}{ 1 +\frac{\gamma-1}{2} {\rm Ma}_{\infty}^2  }.
$$
Figure \ref{fig:bumpflow_vis}(b) shows an horizontal slice of the Mach number at $x_2=0.6$ for three parameters ${\mu}_{\rm min} = [0.75,1.7]$, 
 $\bar{{\mu}} = [0.775,1.75]$
 ${\mu}_{\rm max} = [0.8,1.8]$;
 Figures \ref{fig:bumpflow_vis}(c) and (d)
 show the 
contour lines of the Mach number for 
 ${\mu}_{\rm min}$ and 
 ${\mu}_{\rm max}$: the   red dots  in the Figures denote 
salient points of the flow for $\mu=\mu_{\rm min}$ and  
 are intended to simplify the comparisons between the two flows. 
 
We resort to a DG discretization based on artificial viscosity. We use the local Lax-Friedrichs flux for the advection term, and the BR2 flux (cf.  \cite{bassi1997high}) for the diffusion term. We consider the piecewise-constant viscosity
\begin{equation}
\label{eq:artificial_viscosity}
\left(
\nu( {U})
\right)_k
\,=\,
c_{\rm visc} \, 
\left(\frac{h_k}{ \texttt{p} } \right)^2    \frac{1}{ |\texttt{D}_k| } \int_{\texttt{D}_k} | \nabla \cdot u  | dx 
\end{equation}
where $h_k = \sqrt{  |\texttt{D}_k|} $ is the characteristic size of the $k$-th element of the mesh and $c_{\rm visc} >0$ is a constant set equal to $c_{\rm visc} = 10$ in the numerical simulations. Note that \eqref{eq:artificial_viscosity} is an example of dilation-based model for the viscosity: we refer to the recent review \cite{yu2020study} for   alternative viscosity models and for extensive comparisons.
 
To estimate the hf solution ${U}_{\mu}^{\rm hf} \in \mathcal{X}_{\rm hf}$, we resort to the pseudo-time continuation strategy proposed in \cite{bassi2010very}.
More in detail, if we denote by 
${\mathbf{R}}_{\mu}: \mathcal{X}_{\rm hf} \to \mathbb{R}^{N_{\rm hf}}$ and
by
${\mathbf{J}}_{\mu}: \mathcal{X}_{\rm hf} \to \mathbb{R}^{N_{\rm hf}, N_{\rm hf}}$
 the hf residual and the hf Jacobian and by 
$\mathbf{M} \in \mathbb{R}^{N_{\rm hf}, N_{\rm hf}}$ the mass matrix, we consider the iterative scheme:
\begin{equation}
\label{eq:pseudo_time_integration}
{U}_{\mu}^{\rm hf,k+1} =
{U}_{\mu}^{\rm hf}+
\Delta t_k
\delta {U}_{\mu}^{\rm hf,k+1},
\quad
{\rm with} \;\; 
\left( \mathbf{M} \, + \, 
\Delta t_k {\mathbf{J}}_{\mu}
(  {U}_{\mu}^{\rm hf, k} )
\right) \; \delta {U}_{\mu}^{\rm hf,k+1}
\, = \, 
-
{\mathbf{R}}_{\mu}(  {U}_{\mu}^{\rm hf, k} )
\quad
k=1,2, \ldots,
\end{equation}
where $\Delta t_k$ is chosen adaptively based on the strategy detailed in
 \cite[Chapter 4]{colombo2011agglomeration}. 
Note that \eqref{eq:pseudo_time_integration}  can be interpreted as a  Newton solver with an adaptive relaxation factor.  

We conclude this section by introducing the purely-geometric map used to deform the mesh in absence of a priori information about the solution: in section \ref{sec:registration}, we introduce  a generalization of this map that takes into account the parametric field of interest.
Towards this end, we define $\widehat{\Omega}=(0,1)^2$ and we introduce the parameterized  Gordon-Hall map (cf. \cite{gordon1973construction})  as
\begin{equation}
\label{eq:gordon_hall}
\begin{array}{rl}
{\Psi}_{\mu}(x)
=
&
\displaystyle{
(1-x_2) c_{\rm btm,\mu}(x_1)  + x_2 c_{\rm top}(x_1)
+(1-x_1) c_{\rm left}(x_2)  + x_1c_{\rm right}(x_2)
} 
\\[3mm]
&
\displaystyle{
-
\left(
(1-x_1)(1-x_2) c_{\rm btm,\mu}(0)
+ 
x_1 x_2  c_{\rm top}(1)
+ x_1(1-x_2) c_{\rm btm,\mu}(1)
+
(1-x_1) x_2 c_{\rm top}(0)
\right),
}
\\
\end{array}
\end{equation}
where $c_{\rm btm}, c_{\rm top},c_{\rm left},c_{\rm right}$ are parameterizations of the bottom, top, left and right boundaries of the domain, respectively. Note that $c_{\rm btm}$ depends on the parameter $\mu$ through the  angle $\alpha$  (cf. Figure \ref{fig:bumpflow_vis}(a)): we build $c_{\rm btm}$   so that the jump discontinuities of its derivative $c_{\rm btm}'$ --- which correspond to the extrema of the bump --- are located at $x_1=0.2$ and $x_1=0.6$ for all parameters. We further define the inverse map ${\Lambda}_{\mu}= {\Psi}_{\mu}^{-1}: \Omega_{\mu} \to \widehat{\Omega}$. We have now the elements to introduce
 the parametric mapping $\Phi^{{\rm geo}} $ such that
\begin{equation}
\label{eq:map_gen_apriori}
\Phi_{\mu}^{{\rm geo}} 
\, = \,
\Psi_{\mu} \circ   \Lambda_{\bar{\mu}},
\end{equation}
where $\bar{\mu}$ is the centroid of $\mathcal{P}$.
Given the mesh $\mathcal{T}_{\rm hf}$, we compute the reference points
$\{ x_j^{\rm hf,ref} = \Lambda_{\bar{\mu}}( x_j^{\rm hf}   )\}_{j=1}^{N_{\rm hf,v}}$; then, for any new value of the parameter, we compute the deformed points of the mesh using the identity
$\Phi_{\mu}^{{\rm geo}} ( x_j^{\rm hf}) = \Psi_{\mu}( x_j^{\rm hf,ref})$ for $j=1,\ldots,N_{\rm hf,v}$.

\section{Methodology}
\label{sec:methods}

In this section, we present the  methodology through the vehicle of the model problem introduced in section \ref{sec:model_problem}. In section \ref{sec:registration}, we present the registration procedure, while in section
\ref{sec:ROM}, we discuss in detail the projection-based MOR  scheme. Finally,  in section \ref{sec:multifidelity}, we illustrate the multifidelity approach to reduce offline costs. We state upfront that the two building blocks of our formulation, registration and LSPG formulation in parameterized geometries, have been extensively discussed in 
\cite{taddei2021registration} and \cite{taddei2020discretize}. 

\subsection{Registration}
\label{sec:registration}

The registration procedure takes as input a mesh
$\mathcal{T}_{\rm hf}$  of $\Omega$, a set of snapshots $\{ (  \mu^k, U^k =U_{\mu^k}^{\rm hf}  ) \}_{k=1}^{n_{\rm train}}$, and returns  a parameterized mapping $\Phi: \Omega \times \mathcal{P} \to \mathbb{R}^2$,
$$
\Phi = \texttt{param{\_}registration} \left(
\mathcal{T}_{\rm hf}, \; 
\{ (  \mu^k, U^k =U_{\mu^k}^{\rm hf}  ) \}_{k=1}^{n_{\rm train}}
\right).
$$
In the remainder of this section, we illustrate the  key  features of the procedure and we provide several comments.

\subsubsection{Spectral maps}

The first step of our registration procedure consists in introducing a class of approximation maps.
 Following 
\cite{taddei2021registration},
we consider mappings of the form
\begin{subequations}
\label{eq:map_gen}
\begin{equation}
\label{eq:map_gen_a}
\texttt{N}(\mathbf{a}; \mu   )
\, = \,
\Psi_{\mu} \circ \widetilde{\Phi}  \circ  \Lambda_{\bar{\mu}},
\quad
\widetilde{\Phi} = \texttt{id} + \varphi,
\;\;
\varphi = \sum_{m=1}^M (\mathbf{a})_m \varphi_m.
\end{equation}
Note that 
$\texttt{N}$ generalizes the map \eqref{eq:map_gen_apriori} in the sense that 
$\texttt{N}(\mathbf{0}; \mu ) = \Phi_{\mu}^{\rm geo}$. 
Here, $\bar{\mu}$ is the centroid of $\mathcal{P}$ and
$\varphi_1,\ldots,\varphi_M$ belong to the polynomial space
\begin{equation}
\label{eq:map_gen_b}
\mathcal{W}_{{\rm hf}} =
\left\{
\varphi \in [\mathbb{Q}_J]^2: \;
\varphi \cdot \widehat{n} |_{\partial \widehat{\Omega}} = 0,
\;\;
\varphi(s,0) =0, \; s\in \{0.2,0.6\}
\right\},
\end{equation}
where $\mathbb{Q}_J$ denotes the space of tensorized polynomials of degree at most $J$ in each variable, $\widehat{n}$ is the outward normal to $\widehat{\Omega}$. In the numerical tests, we consider $J=15$.
Note that the second  condition  in \eqref{eq:map_gen_b}  ensures that 
jump discontinuities of $\nabla \texttt{N}(\mathbf{a}; \mu   )$ are located in $[-0.5,0], [0.5,0]$ for all $\mathbf{a}\in \mathbb{R}^M$ and $\mu\in \mathcal{P}$.
\end{subequations}
We equip the mapping space $\mathcal{W}_{{\rm hf}}$ with the $H^2$ norm, 
 \begin{equation}
 \label{eq:norm_mapping_space}
\| \varphi \|_{H^2(\widehat{\Omega})}^2 : =
 \int_{\widehat{\Omega}} \; \left( 
 \sum_{i,j,k=1}^2
( \partial_{j,k} \varphi_i  )^2 \; + \;
 \sum_{i=1}^2
 \varphi_i ^2 \right) \, d x.
 \end{equation}

Exploiting the analysis in \cite{taddei2020registration,taddei2021registration}, we find that
$\mathbf{N}(\mathbf{a}; \mu)$ is a bijection from $\Omega$ to $\Omega_{\mu}$ for  all $\mathbf{a}$ in the set 
\begin{subequations}
\label{eq:calA_bj}
\begin{equation}
\label{eq:calA_bj_true}
\mathcal{A}_{\rm bj} := \left\{
\mathbf{a} \in \mathbb{R}^M \, : \,
\inf_{x \in \widehat{\Omega}} \, 
\widehat{g}({x}; \mathbf{a}) > 0
\right\},
\quad
\widehat{g}(\cdot ;  \mathbf{a})  : = {\rm det} \nabla \widetilde{\Phi}(\mathbf{a}).
\end{equation}
The set $\mathcal{A}_{\rm bj}$ is difficult to deal with numerically: as a result, we define
$\mathcal{A}_{\rm bj}' := \left\{
\mathbf{a} \in \mathbb{R}^M \, : \,
\mathfrak{C}(  \mathbf{a} ) \leq 0 
\right\}$ such that 
\begin{equation}
\label{eq:bijectivity_constraint}
\mathfrak{C}(  \mathbf{a} ) \; : = \;
\int_{\widehat{\Omega}} \;
{\rm exp} \left(
\frac{\epsilon  -   \widehat{g}({x}; \mathbf{a}) }{  C_{\rm exp}  }
\right)
\,+ \,
{\rm exp} \left(
\frac{ \widehat{g}({x}; \mathbf{a}) - 1/\epsilon   }{  C_{\rm exp}  }
\right)
\; dx
-\delta,
\end{equation}
where $\epsilon,C_{\rm exp}, \delta$ are positive constants that will be specified in the next section.
Provided that ${\rm exp} \big(\frac{\epsilon  }{  C_{\rm exp}  }
\big)$ is sufficiently large, we find that there exists a constant $C>0$ such that (see \cite[section 2.2]{taddei2020registration}):
\end{subequations}
\begin{equation}
\label{eq:bijectivity_result}
\mathcal{A}_{\rm bj}  \subset 
\mathcal{A}_{\rm bj}' \cap \{ \mathbf{a} \,  : \, \sup_{x\in \widehat{\Omega}}  \|\nabla  \widehat{g}({x}; \mathbf{a} \|_2 \leq C \}.
\end{equation}
 The discussion above motivates  the combination of the constraint
 $\mathfrak{C}(  \mathbf{a} ) \leq 0 $ with a (strong or weak) control of the second-order derivatives of the mapping.  
We refer to $\mathfrak{C}(  \mathbf{a} ) \leq 0 $ as to the bijectivity constraint.

\subsubsection{Optimization-based registration}

Given $\mu \in \mathcal{P}$, we denote by $s_{\mu} \in L^2(\widehat{\Omega})$ a target sensor that depends on the solution $U_{\mu}$, and we introduce the  $N$-dimensional template space $\mathcal{S}_N \subset L^2(\widehat{\Omega})$. We further denote by $\mathcal{W}_M \subset \mathcal{W}_{\rm hf}$ an $M$-dimensional mapping space and by $W_M: \mathbb{R}^M \to \mathcal{W}_M$ an isometry such that
$\|W_M \mathbf{a}  \|_{H^2(\widehat{\Omega})} = \| \mathbf{a} \|_2$ for all $\mathbf{a} \in \mathbb{R}^M$.
We discuss the construction of 
$\mathcal{S}_N,\mathcal{W}_M$ and the sensor $s_{\mu}$ in the next sections.

We can then introduce the optimization statement that is used to identify the mapping coefficients for a given $\mu\in \mathcal{P}$:
\begin{subequations}
\label{eq:optimization_statement}
\begin{equation}
\begin{array}{l}
\displaystyle{
\min_{\mathbf{a} \in \mathbb{R}^M}
\mathfrak{f}( \mathbf{a};   s_{\mu}, \mathcal{S}_N, W_M  ) 
\; + \;
\xi | W_M \mathbf{a}  |_{H^2(\widehat{\Omega})}^2
\; + \;
\xi_{\rm msh} \mathfrak{R}_{\rm msh} (\mathbf{a};  \mu );
} \\[3mm]
\displaystyle{
{\rm subject \; to} \;\;
\mathfrak{C}(\mathbf{a}) \leq 0,
}
\\
\end{array}
\end{equation}
where
$ | \varphi  |_{H^2(\widehat{\Omega})}^2 =
\int_{\widehat{\Omega}} \;
 \sum_{i,j,k=1}^2
( \partial_{j,k} \varphi_i  )^2 \; + \;
 \sum_{i=1}^2
 \varphi_i ^2 \, d x
$ is the $H^2$ seminorm.
Here, the proximity measure $\mathfrak{f}$ measures the projection error associated with the mapped target $s_{\mu}$ with respect to the template space $\mathcal{S}_N$,
\begin{equation}
\label{eq:proximity_measure}
\mathfrak{f}( \mathbf{a};   s_{\mu}, \mathcal{S}_N, W_M  ) \; := \;
\min_{{\psi} \in \mathcal{S}_N }
\int_{\widehat{\Omega}} \;
\left(
s_{\mu} \circ \tilde{\Phi}(\cdot; \mathbf{a})
- \psi
\right)^2  \; d {x},
\quad
\tilde{\Phi} = \texttt{id} + W_M \mathbf{a}.
\end{equation}
The contribution  $\xi | W_M \mathbf{a}  |_{H^2(\widehat{\Omega})}^2$ is a regularization term that is intended to  control the norm of the mapping Hessian and, in particular, the gradient of the Jacobian determinant $\nabla \hat{g}(\cdot; \mathbf{a})$:
recalling \eqref{eq:bijectivity_result},  the latter is important to enforce bijectivity.
The   term $\mathfrak{R}_{\rm msh}$ penalizes excessive distortions of the mesh and ultimately preserves the discrete bijectivity
(cf. section \ref{sec:notation}):
\begin{equation}
\label{eq:registration_statement_Rmsh}
\mathfrak{R}_{\rm msh} (\mathbf{a}; \mu ) 
=
\sum_{k=1}^{N_{\rm e}} \; |\texttt{D}_k |
{\rm exp} \left(
\mathfrak{f}_{\rm msh,k}\left( 
 \texttt{N} ( \mathbf{a}; \mu)
\right)
\,-\,
\mathfrak{f}_{\rm msh,max}
\right),
\end{equation}
where $\mathfrak{f}_{\rm msh,max}>0$ is a given positive constant and 
\begin{equation}
\label{eq:local_mesh_distortion}
\mathfrak{f}_{\rm msh,k}(\Phi)
: = 
\frac{1}{2}
\frac{\|  \nabla {\Psi}_{k,\Phi}^{\rm hf,1}   \|_{\rm F}^2 }{
|  {\rm det} (  \nabla {\Psi}_{k,\Phi}^{\rm hf,1}   )  |},
\;k=1,\ldots,N_{\rm e},
\end{equation}
$\| \cdot \|_{\rm F}$ is the Frobenius norm and 
${\Psi}_{k,\Phi}^{\rm hf,1} $ is the elemental mapping  associated with  the mapped mesh and a \texttt{p}=1 discretization. We observe that the indicator  \eqref{eq:local_mesh_distortion} is widely used for high-order mesh generation, and has also been considered in  \cite{zahr2020implicit} to prevent mesh degradation, in the DG framework.
Finally, $\mathfrak{C}$ is the bijectivity constraint in \eqref{eq:bijectivity_constraint}. 
\end{subequations}

We observe that the optimization statement depends on several parameters: here, we set
$$
\epsilon=0.1, \;\;
C_{\rm exp} =0.025 \epsilon,  \;\;
\delta = 1, \; \;
\mathfrak{f}_{\rm msh,max}=10, \;\;
\xi = 10^{-3}, \;\;
\xi_{\rm msh} = 10^{-6}.
$$
Since the optimization statement \eqref{eq:optimization_statement} is highly nonlinear and non-convex, the choice of the initial condition  is of paramount importance: here, we exploit the strategy described in \cite[section 3.1.2]{taddei2020registration} to initialize the optimizer; furthermore, we resort to the Matlab function \texttt{fmincon} \cite{MATLAB:2020b}, which relies on an interior penalty algorithm to find local minima of \eqref{eq:optimization_statement}. In our implementation, we provide gradients of the objective function and we rely on a structured mesh on $\widehat{\Omega}$ to speed up evaluations of the sensor  and its gradient at deformed quadrature points,  at each iteration of the optimization algorithm.

\begin{remark}
\label{remark:choice_xi}
In our experience, the choice of $\xi$ is of paramount importance for performance. Small values of $\xi$ lead to lower values of the proximity measure at the price of more irregular mappings (i.e., larger values of $|W_M \mathbf{a}  |_{H^2}$). 
We empirically observe that the latter reduces the generalization properties of the regression algorithm (cf. section \ref{sec:generalization}) used to define the parameterized mapping; in terms of reconstruction performance, we also find that the mapping process introduces small-amplitude smaller spatial scale distortions that ultimately control convergence of the ROM (cf. \cite[Figure 5]{taddei2020registration}) and become more and more noticeable as $\xi$ decreases.
\end{remark}

\subsubsection{Parametric registration}

Given snapshots of the sensor $s$, $\{  (\mu^k, s_{\mu^k} ) \}_{k=1}^{n_{\rm train}}$, we propose to iteratively build the template space $\mathcal{S}_N$, $\mathcal{W}_M$ through the Greedy procedure provided in Algorithm \ref{alg:registration}.
The algorithm takes as input (i) the sensors associated with the snapshot set, (ii)  the initial template $\mathcal{S}_{N_0}$, and 
(iii) the mesh $\mathcal{T}_{\rm hf}$,
 and returns
 (i) the final template space $\mathcal{S}_{N}$,
 (ii) the isometry $W_M$ associated with the mapping space, and 
 (iii) the mapping coefficients $\{  \mathbf{a}^k  \}_k$.
  To clarify the procedure, we introduce  notation
$$
\left[ \mathbf{a}^{\star}, \mathfrak{f}_{N,M}^{\star}  \right]
\, = \,
\texttt{registration} \left(
s,  \mathcal{S}_N, W_M, \mathcal{T}_{\rm hf}, \mu
\right)
$$
to refer to the function that takes as input the target sensor $s$, the template space $\mathcal{S}_N$, the isometry $W_M: \mathbb{R}^M \to \mathcal{W}_M$ associated with the mapping space, the mesh 
$ \mathcal{T}_{\rm hf}$ of $\Omega$ and the parameter $\mu \in \mathcal{P}$ and returns a solution to \eqref{eq:optimization_statement} and the value of the proximity measure
$\mathfrak{f}_{N,M}^{\star} = \mathfrak{f}( \mathbf{a}^{\star}, s, \mathcal{S}_N, W_M   )$. Furthermore, we introduce the POD function that takes as  input a set of mapping coefficients and returns the reduced isometry and the projected mapping coefficients
$$
[  W_M , \; \{  \mathbf{a}^k  \}_k ]  =
\texttt{POD} \left( 
\{ W_{\tilde{M}} \tilde{\mathbf{a}}^{k}  \}_{k=1}^{n_{\rm train}}, 
tol_{\rm pod}  , \| \cdot \|_{H^2(\widehat{\Omega})} \right),
$$
where $M$ is chosen according to the eigenvalues 
$\{  \lambda_m  \}_m$
of the Gramian matrix $\mathbf{C} \in  \mathbb{R}^{n_{\rm train}, n_{\rm train}}$ such that
$\mathbf{C}_{k,k'} =  \tilde{\mathbf{a}}^{k} \cdot  \tilde{\mathbf{a}}^{k'} $, 
\begin{equation}
\label{eq:POD_cardinality_selection}
M := \min \left\{
M': \, \sum_{m=1}^{M'} \lambda_m \geq  \left(1 - tol_{\rm pod} \right) 
\sum_{i=1}^{n_{\rm train}} \lambda_i
\right\}.
\end{equation}

We observe that our approach depends on several hyper-parameters. In our tests, we set  
 $\mathcal{S}_{N_0=1} = {\rm span}\{ s_{\bar{\mu}}  \}$,
 where $\bar{\mu}$ is the centroid of $\mathcal{P}$; furthermore, we set $N_{\rm max}=5$,  $tol_{\rm pod}=10^{-3}$ and  
$\texttt{tol}=10^{-4}$.

\begin{algorithm}[H]                      
\caption{Registration algorithm}     
\label{alg:registration}     

 \small
\begin{flushleft}
\emph{Inputs:}  $\{ (\mu^k,  s^k = s_{\mu^k}) \}_{k=1}^{n_{\rm train}} \subset  \mathcal{P} \times L^2(\widehat{\Omega})$ snapshot set, 
$\mathcal{S}_{N_0} = {\rm span} \{ {\psi}_n \}_{n=1}^{N_0}$ initial template space;
$\mathcal{T}_{\rm hf}$ mesh.
\smallskip

\emph{Outputs:} 
${\mathcal{S}}_N = {\rm span} \{ {\psi}_n  \}_{n=1}^N$ template space, 
$W_M: \mathbb{R}^M \to \mathcal{W}_M$
mapping isometry,
$\{  \mathbf{a}^k  \}_k$ mapping coefficients.
\end{flushleft}                      

 \normalsize 

\begin{algorithmic}[1]
\State
Initialization: 
$\mathcal{S}_{N=N_0} = \mathcal{S}_{N_0}$,
$\mathcal{W}_M =\mathcal{W}_{\rm hf}$.
\vspace{3pt}

\For {$N=N_0, \ldots, N_{\rm max}-1$ }

\State
$
\left[ \mathbf{a}^{\star,k}, \mathfrak{f}_{N,M}^{\star,k}  \right]
\, = \,
\texttt{registration} \left(
s^{k}, \mathcal{S}_N, W_M, \mathcal{T}_{\rm hf}, \mu
\right)$ for $k=1,\ldots,n_{\rm train}$.
\vspace{3pt}

\State
$[  W_M , \; \{  \mathbf{a}^k  \}_k ]  =
\texttt{POD} \left( 
\{ W_M \mathbf{a}^{\star,k}  \}_{k=1}^{n_{\rm train}}, 
tol_{\rm pod}  , \| \cdot \|_{H^2(\widehat{\Omega})}  \right)$
\vspace{3pt}

\If{   $\max_k   \mathfrak{f}_{N,M}^{\star,k}   < \texttt{tol}$}, \texttt{break}

\Else

\State
$\mathcal{S}_{N+1} = \mathcal{S}_{N} \cup {\rm span} 
\{ s^{k^{\star}}     \circ  
(\texttt{id} + W_M \mathbf{a}^{k^{\star}}     )
 \}$ with 
$k^{\star} = {\rm arg} \max_k \mathfrak{f}_{N,M}^{\star,k}$.
\EndIf
\EndFor
\end{algorithmic}
\end{algorithm}

\subsubsection{Choice of the registration sensor}

The sensor $s: \mathcal{P} \to L^2(\widehat{\Omega})$ should be designed to capture relevant features of the solution field that are important to track through registration; furthermore, it should be sufficiently smooth to allow efficient applications of gradient-based  optimization methods. Given the FE field in the deformed mesh   $(\Phi_{\mu}^{\rm geo}( \mathcal{T}_{\rm hf} ), \mathbf{U}_{\mu}^{\rm hf})$, we compute the Mach number ${\rm Ma}_{\mu}^{\rm hf}$ (see 
\eqref{eq:a_mach}) in the nodes of the mesh; then, we define the sensor as the solution to the following smoothing problem:
\begin{equation}
\label{eq:sensor_smoothing}
s_{\mu} := {\rm arg} \min_{s \in H^1(\widehat{\Omega})} \;
\xi_{\rm s} \| \nabla s  \|_{L^2(\widehat{\Omega})}^2
\, + \,
\sum_{i=1}^{N_{\rm hf,v}} 
\left(
s(  {x}_j^{\rm hf, ref}  )
-
\left( \boldsymbol{\rm Ma}_{\mu}^{\rm hf} \right)_j
\right)^2.
\end{equation}
The regularization term associated with the hyper-parameter  $\xi_{\rm s}>0$  is needed due to the fact that $s_{\mu} $ is  defined over  a structured\footnote{As explained in \cite{taddei2021registration}, the use of structured meshes for the sensor is crucial to speed up the evaluation of the objective function of \eqref{eq:optimization_statement}.} mesh over $\widehat{\Omega}$  that is not related to the mesh $\mathcal{T}_{\rm hf}$ used for FE calculations. In all our tests, we consider 
$\xi_{\rm s} =10^{-4}$.
We refer to \cite[section 3.3]{taddei2021registration} for an alternative strategy for the construction of the sensor.

\begin{figure}[h!]
\centering

\subfloat[$\mu = {\mu}_{\rm min}$]{
\includegraphics[width=0.48\textwidth]
{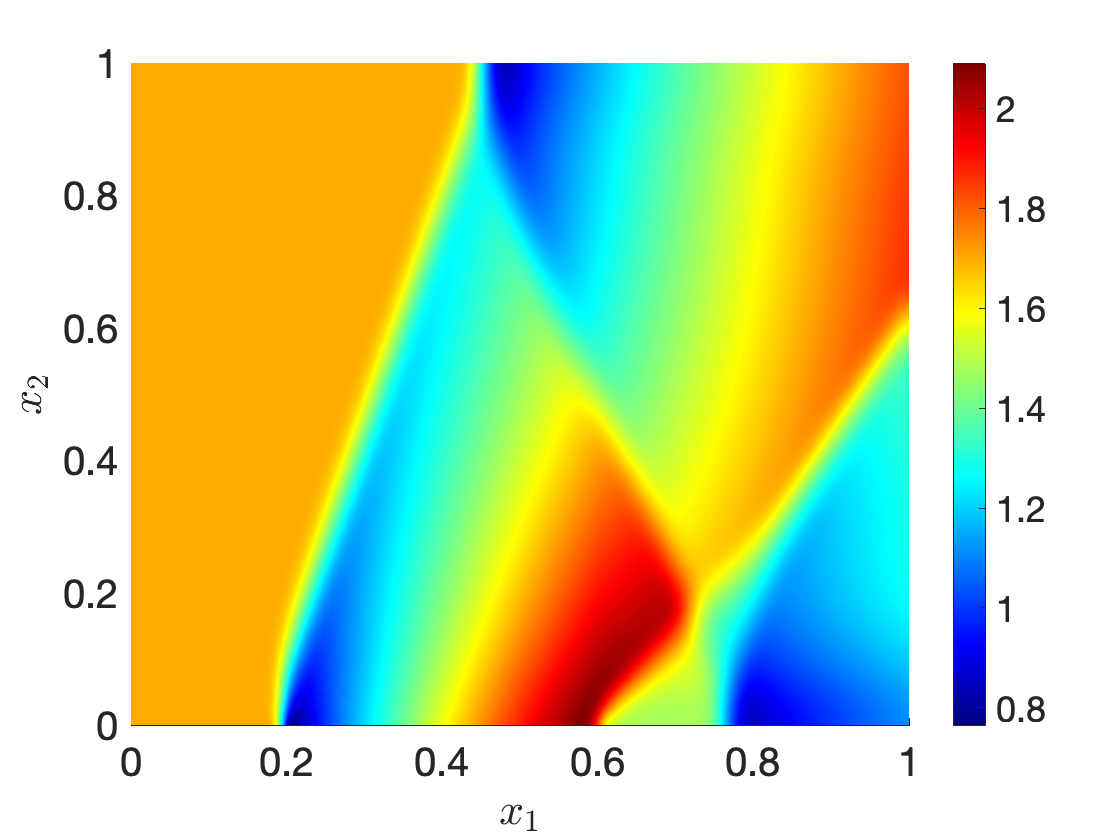}
}
~~~
\subfloat[$\mu = {\mu}_{\rm max}$ ]{
\includegraphics[width=0.48\textwidth]
{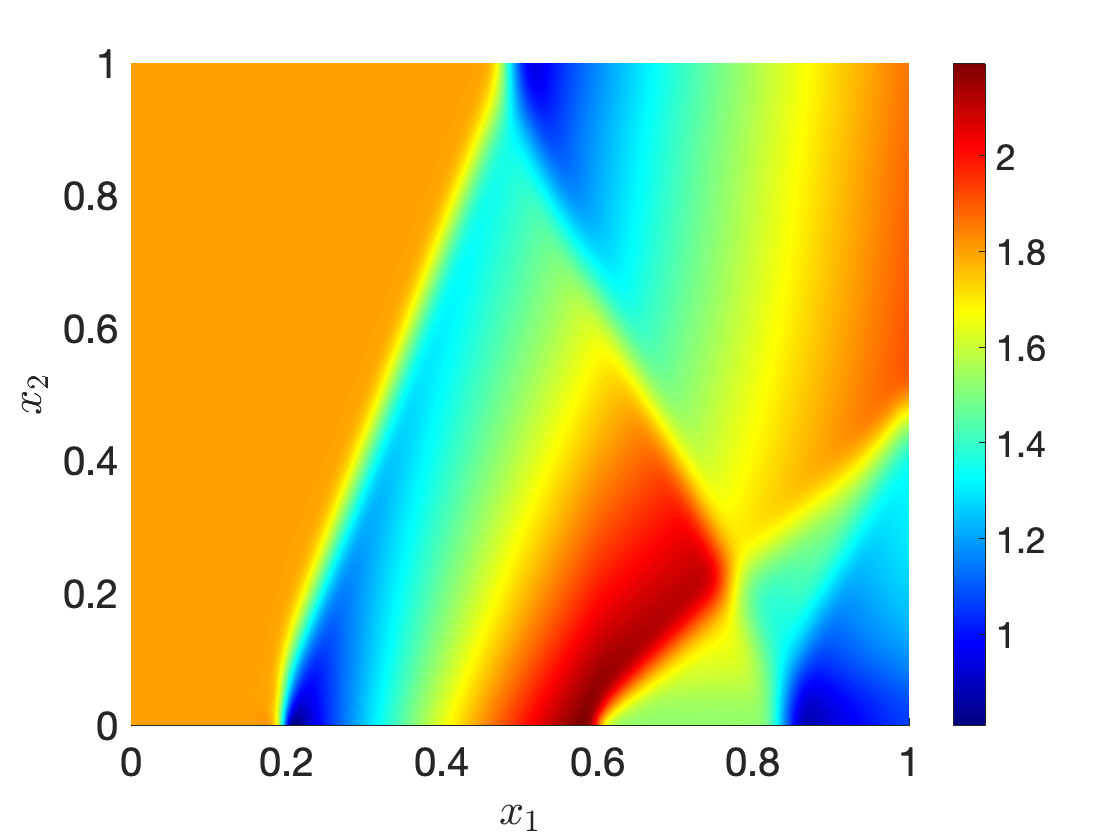}
}

\caption{registration sensor for two values of the parameter, ${\mu}_{\rm min} = [0.75,1.7]$, 
 ${\mu}_{\rm max} = [0.8,1.8]$.
}
\label{fig:sensors}
\end{figure}

We observe that the choice of the Mach number to define the registration sensor is coherent with the choice made in \cite{persson2006sub} to define the highest-modal decay artificial viscosity. Other choices are possible: in particular, using the fluid density $\rho$ in \eqref{eq:sensor_smoothing} as opposed to   ${\rm Ma}$, we obtain
 similar results. Figure \ref{fig:sensors} shows the behavior of the registration sensors for the two values of the parameter in Figure \ref{fig:bumpflow_vis}.

\subsubsection{Generalization}
\label{sec:generalization}
Given the dataset $\{ (\mu^k, \mathbf{a}^k)  \}_{k=1}^{n_{\rm train}}$ as provided by Algorithm \ref{alg:registration}, 
we resort to a multi-target regression algorithm to learn a regressor $\mu \mapsto \widehat{\mathbf{a}}_{\mu}$ and ultimately the  parametric mapping
\begin{equation}
\label{eq:parametric_mapping_Phi}
 {\Phi} : \Omega \times \mathcal{P}  \to \mathbb{R}^2,
\quad
 {\Phi}_{\mu} : =
\texttt{N}(W_M \widehat{\mathbf{a}}_{\mu} ; \mu).
\end{equation}
We here resort to radial basis function (RBF, \cite{wendland2004scattered}) approximation: other regression algorithms could also be considered. 
To avoid overfitting, we retain exclusively modes for which the out-of-sample R-squared is above a given threshold (here, $R_{\rm min}=0.75$): we refer to 
\cite{taddei2020registration,taddei2021space} for further details.

We observe that purely data-driven regression techniques  do not enforce bijectivity for out-of-sample parameters: in practice, we should thus consider sufficiently large training sets in Algorithm \ref{alg:registration}.
This is a major limitation of the present approach that motivates the multifidelity proposal discussed in section \ref{sec:multifidelity}.

\subsection{Projection-based reduced-order model}
\label{sec:ROM}

To clarify the formulation and also provide insights  into the implementation, we introduce a number of definitions and further notation.
Given the FE vector $\mathbf{w}\in \mathbb{R}^{N_{\rm hf}}$, we define the elemental restriction operators $\mathbf{E}_k: \mathbb{R}^{N_{\rm hf}} \to \mathbb{R}^{n_{\rm lp} \cdot D}$ 
such that $\mathbf{E}_k \mathbf{w}$ contains the values of the FE field in the nodes of the $k$-th element for 
$k=1,\ldots,N_{\rm e}$; the elemental restriction operators  
 $\mathbf{E}_k^{\rm ext}: \mathbb{R}^{N_{\rm hf}} \to \mathbb{R}^{n_{\rm lp} \cdot D, 3}$  such that $\mathbf{E}_k^{\rm ext} \mathbf{w}$ contains the values of the FE field in the nodes of  the neighbors of the $k$-th element, for $k=1,\ldots,N_{\rm e}$. We further introduce the set of mesh nodes associated with the $k$-th element and its neighbors:
 $\texttt{X}_k^{\rm hf} = \{ x_{\texttt{T}_{i,k}}^{\rm hf}  \}_{i=1}^{n_{\rm lp}}$ and
  $\texttt{X}_{k,\rm ext}^{\rm hf} = \{ x_{\texttt{T}_{i,k'}}^{\rm hf} :
i=1,\ldots,n_{\rm lp}, \; \overline{\texttt{D}}_k   \cap \overline{\texttt{D}}_{k'} \neq \emptyset \}$; given the mapping $\Phi$, we define
$\Phi_{\mu}( \texttt{X}_k^{\rm hf} )= \{ \Phi_{\mu}(  x_{\texttt{T}_{i,k}}^{\rm hf} )  \}_{i=1}^{n_{\rm lp}}$ and
  $\Phi_{\mu}( \texttt{X}_{k,\rm ext}^{\rm hf} )  = \{ \Phi_{\mu}(  x_{\texttt{T}_{i,k'}}^{\rm hf} ) :
i=1,\ldots,n_{\rm lp}, \; \overline{\texttt{D}}_k   \cap \overline{\texttt{D}}_{k'} \neq \emptyset \}$.

We have now the elements to introduce the DG residual associated with \eqref{eq:conservation_law}:
\begin{subequations}
\label{eq:DG_residual}
\begin{equation}
\label{eq:DG_residual_a}
R_{\mu}^{\rm hf} (\mathbf{U}, \mathbf{V})  \; = \;
\sum_{k=1}^{N_{\rm e}} \;
r_{\mu}^k (  \mathbf{U}, \mathbf{V}  ),
\quad
\forall \; \mathbf{U}, \mathbf{V} \in  \mathbb{R}^{N_{\rm hf}},
\end{equation}
where the local residual $r_{\mu}^k$ corresponds to the contribution to the global residual associated with the $k$-th element of the mesh and depends on the value of the FE fields $U,V$ in the $k$-th element and in its neighbors,
\begin{equation}
\label{eq:DG_residual_b}
r_{\mu}^k (  \mathbf{U}, \mathbf{V}  ) \; = \;
r_{\mu} \left(
 \mathbf{E}_k \mathbf{U}, \,  \mathbf{E}_k \mathbf{V}, \, 
 \mathbf{E}_k^{\rm ext} \mathbf{U}, \,  \mathbf{E}_k^{\rm ext} \mathbf{V}, \; 
 \Phi_{\mu}( \texttt{X}_k^{\rm hf} ), \; 
 \Phi_{\mu}( \texttt{X}_{k,\rm ext}^{\rm hf} ) 
\right),
\quad
k=1,\ldots,N_{\rm e}.
\end{equation}
In the DG literature, schemes in which the primal unknown is only coupled with the unknowns of the adjacent elements are referred to as ``compact'':
the BR2 flux considered in this work is an example of compact treatment of second-order terms for DG formulations (cf. \cite{bassi2005discontinuous}).
Decomposition of the residual as the sum of local elemental contributions is at the foundation of the hf assembling and also of the hyper-reduction procedure discussed below.
We emphasize that the decomposition of the facets' contributions is not unique: 
in order to ensure certain stability and conservation properties for  the hyper-reduced ROM, 
we here consider the energy-stable element-wise decomposition in \cite[section 3.1]{yano2019discontinuous}.
\end{subequations}

Given the reduced-order bases (ROBs) $\mathbf{Z} \in \mathbb{R}^{N_{\rm hf}, N}$ and $\mathbf{Y} \in \mathbb{R}^{N_{\rm hf}, J_{\rm es}}$,  $N\leq J_{\rm es}$, and the trial and test norms $\| \cdot \|$ and $\vertiii{\cdot}$, the EQ-LSPG ROM considered in this work reads as follows: find $\widehat{\mathbf{U}}_{\mu} = \mathbf{Z} \widehat{\boldsymbol{\alpha}}_{\mu}$ to minimize 
\begin{subequations}
\label{eq:LSPG_ROM}
\begin{equation}
\label{eq:LSPG_ROM_a}
\min_{ \boldsymbol{\zeta} \in {\rm col}(\mathbf{Z})    }
\;
\sup_{ \boldsymbol{\eta} \in {\rm col}(\mathbf{Y})    }
\;
\frac{R_{\mu}^{\rm eq}( \boldsymbol{\zeta} , \boldsymbol{\eta}     ) }{
\vertiii{ \boldsymbol{\eta} }  }.
\end{equation}
Here, $R_{\mu}^{\rm eq}$ is the empirical residual defined as 
\begin{equation}
\label{eq:LSPG_ROM_b}
R_{\mu}^{\rm eq} (\mathbf{U}, \mathbf{V})  \; = \;
\sum_{k\in \texttt{I}_{\rm eq}} \; \rho_k^{\rm eq}
r_{\mu}^k (  \mathbf{U}, \mathbf{V}  ),
\quad
\forall \; \mathbf{U}, \mathbf{V} \in  \mathbb{R}^{N_{\rm hf}},
\end{equation}
where $\texttt{I}_{\rm eq}\subset \{1, \ldots,N_{\rm e}\}$ are the indices  of the sampled elements and  
$\boldsymbol{\rho}^{\rm eq} = [\rho_1^{\rm eq},\ldots,\rho_{N_{\rm e}}^{\rm eq}]$
are positive empirical weights to be determined,
$\rho_k^{\rm eq}> 0 \Leftrightarrow k \in \texttt{I}_{\rm eq} $. Provided that the columns $[\boldsymbol{\eta}_1,\ldots,  \boldsymbol{\eta}_{J_{\rm es}} ]$ of $\mathbf{Y}$ are orthonormal with respect to the $\vertiii{\cdot}$ norm, we can rewrite \eqref{eq:LSPG_ROM_a} as
\begin{equation}
\label{eq:LSPG_ROM_c}
 \widehat{\boldsymbol{\alpha}}_{\mu} \in {\rm arg} \min_{\boldsymbol{\alpha} \in \mathbb{R}^N} \; 
 \|  \mathbf{R}_{\mu}^{\rm eq}( \boldsymbol{\alpha} )  \|_2,
\quad
  \mathbf{R}_{\mu}^{\rm eq}( \boldsymbol{\alpha} )  
=\left[
R_{\mu}^{\rm eq}(\mathbf{Z} \boldsymbol{\alpha},  \boldsymbol{\eta}_1), \ldots,
R_{\mu}^{\rm eq}(\mathbf{Z} \boldsymbol{\alpha},  \boldsymbol{\eta}_{J_{\rm es}})
\right]. 
\end{equation}
Note that \eqref{eq:LSPG_ROM_c} is a nonlinear least-squares problem that can be solved using the Gauss-Newton algorithm.
We initialize the iterative procedure using a non-intrusive estimate of the solution coefficients:
if the number of training points is sufficiently large --- such as in the case of POD data compression ---  we use 
RBF regression as in \cite{taddei2021space,taddei2020discretize}; for small training sets --- such as in the first steps of the Greedy algorithm ---  we use nearest-neighbors regression.
Similarly to \cite{taddei2021space}, we resort to a discrete $L^2$ norm for the trial space and to a discrete $H^1$ norm for the test space: we refer to \cite{brunken2019parametrized} for a discussion  on variational  formulations for first-order linear  hyperbolic problems.
\end{subequations}

The MOR formulation \eqref{eq:LSPG_ROM} depends on the choice of the trial and test  ROBs $\mathbf{Z}$ and $\mathbf{Y}$ and on the sparse vector of empirical weights $\boldsymbol{\rho}^{\rm eq}$: we discuss their construction in the remainder of section \ref{sec:ROM}. Before proceeding with the discussion, we remark that we can exploit \eqref{eq:LSPG_ROM_b}   to assemble the reduced residual  
$\mathbf{R}_{\mu}^{\rm eq}$: first, we 
evaluate  $\Phi_{\mu}( \texttt{X}_k^{\rm hf} )$ and   $\Phi_{\mu}( \texttt{X}_{k,\rm ext}^{\rm hf} )$ for all $k\in \texttt{I}_{\rm eq}$; then, we compute the local residuals
$\{ r_{\mu}^k(\mathbf{Z} \boldsymbol{\alpha}, \boldsymbol{\eta}_j)\}_k$ using \eqref{eq:DG_residual_b}; finally, we compute $ \mathbf{R}_{\mu}^{\rm eq}( \boldsymbol{\alpha} )$  by summing over the sampled elements, cf. 
\eqref{eq:LSPG_ROM_b}. Note that, since the residuals $\{ r_{\mu}^k\}_k$ are linear with respect to the test function, we can use standard element-wise residual evaluation routines to compute local contributions to the residual. Furthermore, we observe that computation of the residual $ \mathbf{R}_{\mu}^{\rm eq}( \boldsymbol{\alpha} )$ requires the storage of trial and test ROBs in the sampled elements  and in their neighbors, and is thus independent of the total number of mesh elements. We refer to \cite{taddei2020discretize} for further details.

\begin{remark}
\label{remark:dtm}
We observe that we here resort to a discretize-then-map (DtM, \cite{dal2019hyper,taddei2020discretize,washabaugh2016use})
treatment of parameterized geometries. As discussed in \cite{taddei2020discretize}, the DtM approach --- as opposed to the more standard map-then-discretize (MtD, 
\cite{lassila2014model,rozza2021basic,ballarin2016fast,ballarin2017numerical,rozza2007reduced}) 
) approach ---  in combination with EQ
 allows to reuse hf local integration routines and is thus considerably easier to implement, particularly for nonlinear PDEs.
\end{remark}

\subsubsection{Construction of trial and test spaces}
 
We resort to the standard data compression algorithms POD and weak-Greedy to build the trial ROB $\mathbf{Z}$.
For stability reasons, we ensure that the columns 
$\boldsymbol{\zeta}_1,\ldots,
\boldsymbol{\zeta}_N$ of 
$\mathbf{Z}$ are orthonormal with respect to the $\| \cdot \|$ norm.
 We anticipate that, for the problem considered in this paper, POD leads to superior performance (cf. section \ref{sec:numerics}) in terms of online accuracy; however, POD requires more extensive explorations of the parameter domain and is thus more onerous during the offline stage. 
For this reason, in section \ref{sec:multifidelity},  we resort to the weak-Greedy method in combination with  multi-fidelity training to reduce offline costs.
We refer to the monographies 
\cite{hesthaven2016certified,quarteroni2015reduced}
for extensive discussions on 
 POD and weak-Greedy data compression.

For completeness, we report in Algorithm \ref{alg:weak_greedy} the 
 weak-greedy algorithm as implemented in our code. 
Note that the algorithm  takes as input the mesh $\mathcal{T}_{\rm hf}$ and the mapping $\Phi$  which define the FE mesh for all parameters, and returns the ROB $\mathbf{Z}$ and the ROM for the solution coefficients.
  The residual indicator is presented in section \ref{sec:dual_residual_estimation}.
 The function \texttt{Gram-Schmidt} at Line 4 performs one step of  the Gram Schmidt process to ensure  that the trial ROB is orthonormal with respect to the $\| \cdot \|$ norm. 
  Construction of the ROM at Line 5 involves the construction of the test ROB $\mathbf{Y}$ and the computation of the empirical quadrature rule: these procedures are described below.

\begin{algorithm}[H]                      
\caption{Weak-greedy algorithm. }     
\label{alg:weak_greedy}     

 \small
\begin{flushleft}
\emph{Inputs:}  $\mathcal{P}_{\rm train} : =   \{ \mu^k \}_{k=1}^{n_{\rm train}} $ training parameter set, 
$\Phi:\Omega\times \mathcal{P} \to \mathbb{R}^2$ mapping;
$\mathcal{T}_{\rm hf}$ mesh.
\smallskip

\emph{Outputs:} 
$\mathbf{Z}$ trial ROB;
$\mu \in \mathcal{P} \mapsto \widehat{\boldsymbol{\alpha}}_{\mu}$
ROM for the solution coefficients.

\end{flushleft}                      

 \normalsize

 \textbf{Offline stage}
\medskip

\begin{algorithmic}[1]
\State
Choose 
$\mu^{\star, 1} = \bar{\mu}$.
\vspace{3pt}

\For {$N=1,\ldots,N_{\rm max}$}

\State
Solve the hf problem for $\mu = \mu^{\star,N}$ to obtain $\mathbf{U}^{\star} = \mathbf{U}_{\mu^{\star,N}}$.
\vspace{3pt}

\State
Update the ROB $\mathbf{Z} =\texttt{Gram-Schmidt}(\mathbf{Z},\mathbf{U}^{\star}, \| \cdot  \| ) $.
\vspace{3pt}

\State
Build the ROM
$\mu \in \mathcal{P} \mapsto \widehat{\boldsymbol{\alpha}}_{\mu}$.
\vspace{3pt}

 \For {$k=1,\ldots,n_{\rm train}$}

\State
Estimate the solution using the ROM for $\mu=\mu^k$.
\vspace{3pt}

\State
Compute the error indicator
$\Delta_{\mu^k}:= \mathfrak{R}_{\mu^k}(   \widehat{\boldsymbol{\alpha}}_{\mu^k}   )$ (cf. section \ref{sec:dual_residual_estimation}).
\EndFor

\State
Set $ \mu^{\star,N+1} = {\rm arg} \max_{ \mu\in \mathcal{P}_{\rm train}  } \Delta_{\mu}$. 
\vspace{3pt}

\EndFor
\end{algorithmic}
\end{algorithm}

As  rigorously proven  in \cite[Appendix C]{taddei2021space} for linear inf-sup stable problems, the test ROB $\mathbf{Y}$   should approximate the Riesz representers of the Fr{\'e}chet derivative of the residual at $\mathbf{U}_{\mu}^{\rm hf}$ applied to the elements of the trial ROB for all $\mu \in \mathcal{P}$. Similarly to \cite{taddei2020discretize}, we here resort to the sampling strategy 
 based on POD proposed in \cite{taddei2021space}: 
first, 
given the $\mathcal{Y}$ inner product $(( \cdot, \cdot))$ such that
$\vertiii{\cdot} = \sqrt{(( \cdot, \cdot))}$,
we compute the Riesz representers of the Fr{\'e}chet derivative of the residual at $\mathbf{U}_{\mu}^{\rm hf}$, evaluated for the elements of the $n$-th trial bases $\boldsymbol{\zeta}_n$ and 
for the $k$-th parameter $\mu^k$ in the training set,
$$
((\boldsymbol{\psi}_{k,n}, \mathbf{v}    ))
\, = \,
D R_{\mu}^{\rm hf}[   \mathbf{U}_{\mu}^{\rm hf}    ]
( \boldsymbol{\zeta}_n, \mathbf{v} ),
\quad
\forall \; \mathbf{v} \in  \mathbb{R}^{N_{\rm hf}},
$$
for $n=1,\ldots,N$, $k=1,\ldots,n_{\rm train}$; then, we apply POD for a given tolerance $tol_{\rm test}>0$ to find the test ROB  $\mathbf{Y}$,
$$
[  \mathbf{Y} , \cdot  ]  =
\texttt{POD} \left( 
\{ \boldsymbol{\psi}_{k,n}  \}_{k,n}, 
tol_{\rm test}  , \vertiii{\cdot} \right).
$$ 
The POD tolerance should be sufficiently tight to ensure the well-posedness of the reduced problem:  
in the numerical tests of section \ref{sec:numerics}, we set $tol_{\rm test} =10^{-3}$.

\subsubsection{Empirical quadrature}

As in \cite{taddei2020discretize}, we seek 
$\boldsymbol{\rho}^{\rm eq} \in \mathbb{R}_+^{N_{\rm e}}$ such that
(i) the number of nonzero entries in 
 $\boldsymbol{\rho}^{\rm eq}$,  $\| \boldsymbol{\rho}^{\rm eq}   \|_{\ell^0}$,   is as small as possible;
(ii, \emph{constant function constraint})
 the constant function is approximated correctly in $\Omega$ (i.e., ${\Phi} = \texttt{id}$), 
 \begin{equation}
 \label{eq:constant_function_constraint}
\Big|
\sum_{k=1}^{N_{\rm e}} \rho_k^{\rm eq} | \texttt{D}_k  |
\,-\,
| \Omega | 
 \Big|
 \ll 1; 
 \end{equation}
(iii, \emph{manifold accuracy constraint})
for all $\mu \in \mathcal{P}_{\rm train,eq} = \{  \mu^k \}_{k=1}^{n_{\rm train}+n_{\rm train,eq}}$, the  empirical residual satisfies
\begin{subequations}
\label{eq:accuracy_constraint}
\begin{equation}
\Big \|
   {\mathbf{R}}_{\mu}^{\rm hf}
( \boldsymbol{\alpha}_{\mu}^{\rm train}   )   
\, - \,
   {\mathbf{R}}_{\mu}^{\rm eq}
( \boldsymbol{\alpha}_{\mu}^{\rm train}   )   
 \Big \|_2
 \ll 1.
\end{equation}
where $ {\mathbf{R}}_{\mu}^{\rm hf}$  corresponds to substitute
$\rho_1^{\rm eq} = \ldots = \rho_{N_{\rm e}}^{\rm eq} = 1$ in \eqref{eq:LSPG_ROM_b} and $\boldsymbol{\alpha}_{\mu}^{\rm train}$ satisfies
\begin{equation}
\boldsymbol{\alpha}_{\mu}^{\rm train} = 
\left\{
\begin{array}{ll}
\displaystyle{
\mathbf{Z}^T \mathbf{X}_{\rm hf} \mathbf{U}_{\mu}^{\rm hf} } & {\rm if} \; \mu \in \mathcal{P}_{\rm train} ; \\[3mm]
\displaystyle{ {\rm arg} \min_{\boldsymbol{\alpha} \in \mathbb{R}^N} \; 
 \|  \mathbf{R}_{\mu}^{\rm hf}( \boldsymbol{\alpha} )  \|_2,}
  & {\rm if} \; \mu \notin \mathcal{P}_{\rm train} .\\
\end{array}
\right.
\end{equation}
Here, $\mathbf{X}_{\rm hf}$ is the matrix associated with the $(\cdot, \cdot)$ inner product and 
$\mathcal{P}_{\rm train} = \{  \mu^k \}_{k=1}^{n_{\rm train}}$ is the set of parameters for which the hf solution is available.
When we apply POD to generate the ROM, we set 
$\mathcal{P}_{\rm train}=\mathcal{P}_{\rm train,eq}$; when we apply the weak-Greedy algorithm, we augment $\mathcal{P}_{\rm train}$ with $n_{\rm train,eq}=10$ randomly-selected parameters (see \cite[Algorithm 1]{yano2019discontinuous}): we empirically observe that this choice improves performance of the hyper-reduced ROM, particularly for small values of $n_{\rm train}$. 
We refer to the above-mentioned literature for a thorough motivation of  the previous constraints; in  particular, we refer to 
\cite{chan2020entropy,yano2019discontinuous} for a discussion on the conservation properties of the ROM for conservation laws.
\end{subequations}

It is possible to show  (see, e.g.,  \cite{taddei2021space})   that (i)-(ii)-(iii) lead to a  sparse representation problem of the form
\begin{equation}
\label{eq:sparse_representation}
\min_{  \boldsymbol{\rho} \in \mathbb{R}^{N_{\rm e}} }
\;
\| \boldsymbol{\rho}   \|_{\ell^0},
\quad
{\rm s.t} \quad
\left\{
\begin{array}{l}
\|\mathbf{G} \boldsymbol{\rho} - \mathbf{b}  \|_2 \leq \delta; \\[3mm]
\boldsymbol{\rho} \geq \mathbf{0}; \\
\end{array}
\right.
\end{equation}
for a suitable threshold  $\delta>0$, and for  a suitable  choice of 
$\mathbf{G}, \mathbf{b}$. Following \cite{farhat2015structure}, we here resort to the non-negative  least-squares  method  to find approximate solutions to \eqref{eq:sparse_representation}. In particular, we use the Matlab  function \texttt{lssnonneq}, which takes as input the pair $(\mathbf{G}, \mathbf{b})$ and a tolerance $tol_{\rm eq}>0$ and returns the sparse vector $\boldsymbol{\rho}^{\rm eq}$,
\begin{equation}
\label{eq:lsqnonneg}
[\boldsymbol{\rho}^{\rm eq}] =  \texttt{lsqnonneg} \left( \mathbf{G}, \mathbf{b}, tol_{\rm eq} \right).
\end{equation}
We refer to \cite{chapman2017accelerated} for an efficient implementation of the non-negative least-squares method for large-scale problems.
 
\subsubsection{Dual residual estimation}
\label{sec:dual_residual_estimation}

We here resort to the dual residual error indicator
\begin{equation}
\label{eq:dual_residual_norm}
\mathfrak{R}_{\mu}^{\rm hf}(\boldsymbol{\alpha})
\, : = \,
\sup_{\mathbf{v} \in \mathbb{R}^{N_{\rm hf}}    }
\;
\frac{ R_{\mu}^{\rm hf}(\mathbf{Z} \boldsymbol{\alpha}, \, \mathbf{v} )}{\vertiii{\mathbf{v}}},
\quad
\boldsymbol{\alpha} \in \mathbb{R}^N,
\end{equation}
to drive the weak-Greedy algorithm.
If we denote by  $\mathbf{Y}_{\rm hf}$ the matrix associated with the $\vertiii{\cdot}$ norm, we have that
$$
\mathfrak{R}_{\mu}^{\rm hf}(\boldsymbol{\alpha})
\, : = \,
\sqrt{
\mathbf{R}_{\mu}^{\rm hf}(\mathbf{Z} \boldsymbol{\alpha} )^T  \, 
\mathbf{Y}_{\rm hf}^{-1} \, 
\mathbf{R}_{\mu}^{\rm hf}(\mathbf{Z} \boldsymbol{\alpha} )
},
\quad
\forall \, \mu\in \mathcal{P},  \boldsymbol{\alpha}\in \mathbb{R}^N.
$$
Computation of 
$\mathfrak{R}_{\mu}^{\rm hf}(\boldsymbol{\alpha})$ thus  requires to assemble the hf residual 
$\mathbf{R}_{\mu}(\boldsymbol{\alpha}) \in  \mathbb{R}^{N_{\rm hf}}$ and then solve a linear problem of size $N_{\rm hf}$.  Since the matrix $\mathbf{Y}_{\rm hf}$  is symmetric positive definite and parameter-independent, we use Cholesky factorization to speed up computations 
of the inner loop in Algorithm \ref{alg:weak_greedy}  --- we further use the Matlab function \texttt{symamd} to reduce fill-in.

In the numerical results (cf. Appendix \ref{sec:further_numerics}), we show that 
$\mathfrak{R}_{\mu}^{\rm hf}(\cdot)$ is highly correlated with the relative error. In order to use 
$\mathfrak{R}_{\mu}^{\rm hf}(\cdot)$ during the online stage, we shall perform hyper-reduction: we refer to \cite{taddei2020discretize} for the details.
In our experience, for the value of $n_{\rm train}$  and for the particular hf discretization considered, the cost of the greedy  search in Algorithm \ref{alg:weak_greedy} is negligible compared to the cost of an hf solve; as a result, hyper-reduction does  not seem  needed during the offline stage.

\subsection{Offline/online computational decomposition based on two-fidelity sampling}
\label{sec:multifidelity}

As discussed in section \ref{sec:registration}, the registration procedure relies on a regression algorithm to compute the mapping coefficients
$\widehat{\mathbf{a}}_{\mu}$ 
 for out-of-sample parameters. Since the regression algorithm does not explicitly ensure that bijectivity is satisfied, in practice we should consider sufficiently large training sets $\mathcal{P}_{\rm train}$. 
To address this issue, we propose to use a multi-fidelity approach, which relies on hf solves on  a coarser grid to learn the parametric mapping $\Phi$. Algorithm \ref{alg:offline_online} summarizes the offline/online procedure implemented in our code.

We state below several  remarks.
\begin{itemize}
\item
The snapshots $\{  U_{\mu^k}^{\rm hf,c}  \}_{k=1}^{n_{\rm train,c}}$ are exclusively used to compute the sensors $\{  s_{\mu^k}  \}_{k=1}^{n_{\rm train,c}}$ that are then fed into the registration algorithm: 
we might then employ snapshots from third-party solvers and we might also use different grids for different parameters.
\item
In this work, we propose to build the fine mesh $\mathcal{T}_{\rm hf}$ based on the coarse snapshot $U_{\bar{\mu}}^{\rm hf,c}$; we use here the open source mesh generator proposed in
\cite{persson2004simple} based on a suitable relative size function: we provide details concerning the definition of the size function in Appendix \ref{sec:mesh_generation}. As anticipated in the introduction, we expect that for more challenging problems it might be necessary to adapt the mesh based on multiple snapshots.
\item
Computation of the ROB $\mathbf{Z}$ and of the ROM for the solution coefficients and the online evaluation can be performed using standard pMOR algorithms for linear approximations  in parameterized geometries: we believe that this  represents a valuable feature of the proposed approach that allows its immediate application to a broad class of problems.
\item
Our multi-fidelity procedure does not include any update of the sensors as more accurate simulations become available during Step 5 of the offline stage: as a result, it might lead to poor results if the initial discretization is excessively inaccurate. Development of more sophisticated multi-fidelity techniques is the subject  of ongoing research.
\end{itemize}

\begin{algorithm}[H]                      
\caption{Offline online algorithm. }     
\label{alg:offline_online}     

 \vspace{3pt}
 
 \textbf{Offline stage}
\medskip

\begin{algorithmic}[1]
\State
Generate the snapshots 
$\{  U_{\mu^k}^{\rm hf,c}  \}_{k=1}^{n_{\rm train,c}}$ based on the grid 
$\mathcal{T}_{\rm hf,c}$ and the mapping $\Phi^{\rm geo}$.
\vspace{3pt}

\State
Use the snapshots 
$\{  U_{\mu^k}^{\rm hf,c}  \}_{k=1}^{n_{\rm train,c}}$ to compute the sensors
$\{  s_{\mu^k}  \}_{k=1}^{n_{\rm train,c}}$ using \eqref{eq:sensor_smoothing}.
\vspace{3pt}

\State
Generate the fine mesh $\mathcal{T}_{\rm hf}$.
\vspace{3pt}

\State
Apply  
registration (cf. Algorithm \ref{alg:registration}) based on $\{  s_{\mu^k}  \}_{k=1}^{n_{\rm train,c}}$  and the mesh $\mathcal{T}_{\rm hf}$.
\vspace{3pt}

\State
Generate  the ROB  $\mathbf{Z}$ and the ROM
for the solution coefficients $\mu\in \mathcal{P}  \mapsto \widehat{\boldsymbol{\alpha}}_{\mu} \in \mathbb{R}^N$.
\vspace{3pt}
\end{algorithmic}

 \textbf{Online stage} (for any given $\mu\in \mathcal{P}$)
\medskip

\begin{algorithmic}[1]
\State
Solve the ROM to compute $\widehat{\boldsymbol{\alpha}}_{\mu}$.
\medskip

\State
Compute the deformed mesh $\Phi_{\mu}(\mathcal{T}_{\rm hf})$ and  
$\widehat{\mathbf{U}}_{\mu} = \mathbf{Z} \widehat{\boldsymbol{\alpha}}_{\mu}$.

\end{algorithmic}

\end{algorithm}

\section{Numerical results}
\label{sec:numerics}

We present below extensive numerical investigations for the model problem introduced in section \ref{sec:model_problem}.
Further numerical tests are provided in Appendix \ref{sec:further_numerics}.

\subsection{Test 1: single-fidelity training}

In this first test, we consider performance of our approach without multi-fidelity training.
Towards this end, we consider a \texttt{p}=2 DG FE discretization with $N_{\rm hf} =  197856$ degrees of freedom ($N_{\rm e}=8204$):
the FE mesh is depicted in Figure \ref{fig:test2_meshing}(a).
We consider an equispaced  grid of $11\times 11$
parameters $\mathcal{P}_{\rm train} : = \{ \mu^k\}_{k=1}^{ n_{\rm train}   } \subset \mathcal{P}$
 ($n_{\rm train}=121$); we further consider 
 $n_{\rm test}=10$ randomly-selected parameters for testing. We measure performance of the ROM in terms of the average out-of-sample relative prediction error :
 \begin{equation}
 \label{eq:Eavg}
 E_{\rm avg} 
:=
\frac{1}{n_{\rm test}}
\sum_{\mu \in \mathcal{P}_{\rm test}}
\;
\frac{ \|U_{\mu}^{\rm hf}    - \widehat{U}_{\mu}^{\rm hf}     \|_{L^2(\Omega_{\mu})} }{
\|U_{\mu}^{\rm hf}   \|_{L^2(\Omega_{\mu})} }.
 \end{equation}
 The  mapping $\Phi$ that is obtained applying the registration procedure in Algorithm \ref{alg:registration} consists of three modes ($M=3$):  the R-squared associated with the RBF regressors is above the threshold for all three modes.  

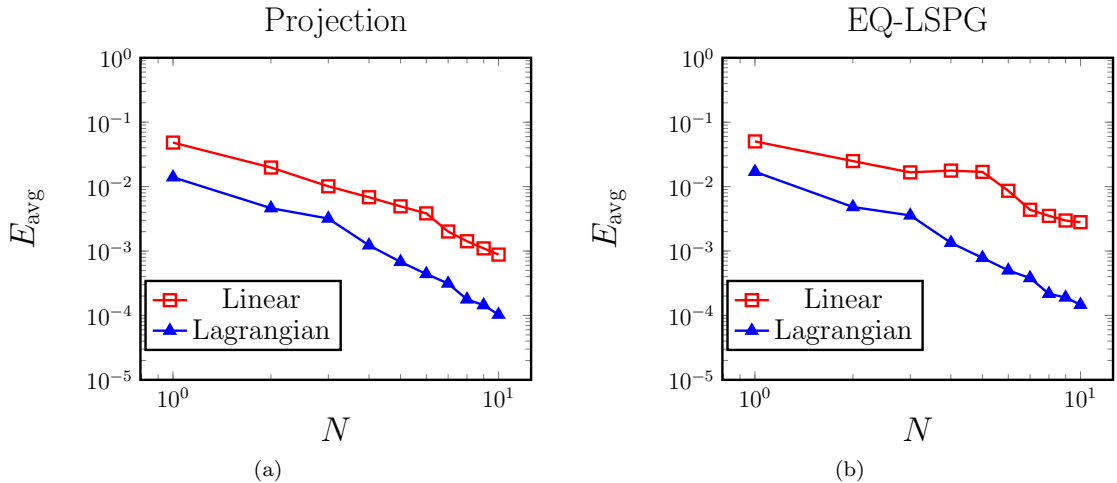
\begin{figure}[h!]
\centering

\subfloat[]{
\begin{tikzpicture}[scale=0.75]
\begin{loglogaxis}[
title = {\LARGE {Projection}},
xlabel = {\LARGE {$N$}},
  ylabel = {\LARGE {$E_{\rm avg}$}},
legend entries = {Linear, Lagrangian},
  line width=1.2pt,
  mark size=3.0pt,
 ymin=0.00001,   ymax=1,
ylabel style = {font=\large,yshift=5ex},
 yticklabel style = {font=\large,xshift=0ex},
xticklabel style = {font=\large,yshift=0ex},
legend style={at={(0.01,0.2)},anchor=west,font=\Large}
  ]
 
  \addplot[color=red,mark=square]  table {data/euler/test1/proj_linear.dat};
  
  \addplot[color=blue,mark=triangle*] table {data/euler/test1/proj_nonlinear.dat};
  
\end{loglogaxis}
\end{tikzpicture}
}
~~~
\subfloat[]{
\begin{tikzpicture}[scale=0.75]
\begin{loglogaxis}[
title = {\LARGE {EQ-LSPG}},
xlabel = {\LARGE {$N$}},
  ylabel = {\LARGE {$E_{\rm avg}$}},
legend entries = {Linear, Lagrangian},
  line width=1.2pt,
  mark size=3.0pt,
 ymin=0.00001,   ymax=1,
ylabel style = {font=\large,yshift=5ex},
 yticklabel style = {font=\large,xshift=0ex},
xticklabel style = {font=\large,yshift=0ex},
legend style={at={(0.01,0.2)},anchor=west,font=\Large}
]
\addplot[color=red,mark=square]  table {data/euler/test1/EQLSPG_linear.dat};
  
  \addplot[color=blue,mark=triangle*] table {data/euler/test1/EQLSPG_nonlinear.dat};
  
\end{loglogaxis}
\end{tikzpicture}
}

\caption{single-fidelity training. Comparison of linear and Lagrangian approaches. Trial ROB $\mathbf{Z}$ is built using POD.
}
\label{fig:test1_pod_training}
\end{figure}

\begin{figure}[h!]
\centering
\subfloat[]{
 \begin{tikzpicture}[scale=0.75]
\begin{semilogyaxis}[
title = {\LARGE {Projection}},
xlabel = {\LARGE {$N$}},
  ylabel = {\LARGE {$E_{\rm avg}$}},
legend entries = {Linear, Lagrangian},
  line width=1.2pt,
  mark size=3.0pt,
xmin=4,   xmax=12,
 ymin=0.00001,   ymax=1,
ylabel style = {font=\large,yshift=5ex},
 yticklabel style = {font=\large,xshift=0ex},
xticklabel style = {font=\large,yshift=0ex},
legend style={at={(0.01,0.2)},anchor=west,font=\Large}
  ]
 
  \addplot[color=red,mark=square]  table {data/euler/test1/proj_linear_greedy.dat};
  
  \addplot[color=blue,mark=triangle*] table {data/euler/test1/proj_nonlinear_greedy.dat};
  
\end{semilogyaxis}
\end{tikzpicture}
}
~~~
\subfloat[]{
 \begin{tikzpicture}[scale=0.75]
\begin{semilogyaxis}[
title = {\LARGE {EQ-LSPG}},
xlabel = {\LARGE {$N$}},
  ylabel = {\LARGE {$E_{\rm avg}$}},
legend entries = {Linear, Lagrangian},
  line width=1.2pt,
  mark size=3.0pt,
xmin=4,   xmax=12,
 ymin=0.00001,   ymax=1,
ylabel style = {font=\large,yshift=5ex},
 yticklabel style = {font=\large,xshift=0ex},
xticklabel style = {font=\large,yshift=0ex},
legend style={at={(0.01,0.2)},anchor=west,font=\Large}
]
\addplot[color=red,mark=square]  table {data/euler/test1/EQLSPG_linear_greedy.dat};
  
  \addplot[color=blue,mark=triangle*] table {data/euler/test1/EQLSPG_nonlinear_greedy.dat};
  
\end{semilogyaxis}
\end{tikzpicture}
}

\caption{single-fidelity training. Comparison of linear and Lagrangian approaches. Trial ROB $\mathbf{Z}$ is built using weak-Greedy.
}
\label{fig:test1_greedy_training}
\end{figure}
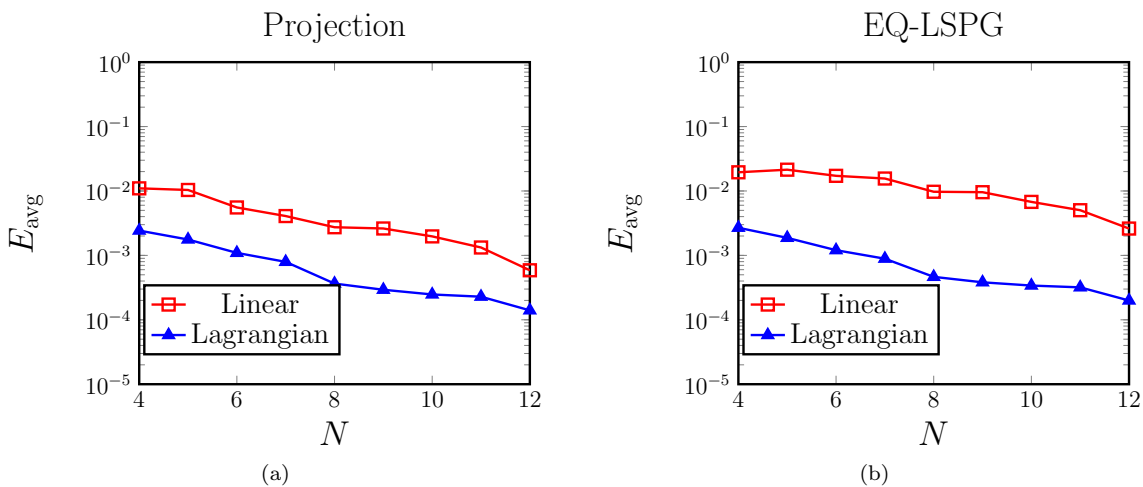

 Figure \ref{fig:test1_pod_training} shows performance of linear and Lagrangian approaches based on POD data compression. Figure \ref{fig:test1_pod_training}(a) shows the projection error, while Figure \ref{fig:test1_pod_training}(b) shows the error associated with the EQ-LSPG ROM introduced in section \ref{sec:ROM}. We observe that registration significantly improves performance for all values of $N$. 
 Figure \ref{fig:test1_greedy_training} replicates the results for the ROM based on weak-Greedy\footnote{We initialize the Greedy procedure with  $N_0=4$ equispaced samples. The Greedy search is performed over the training set of $n_{\rm train}=121$ parameters.} compression: note that also in this case registration significantly improves performance for all values of $N$ considered.
  We further observe that our EQ-LSPG ROM is able to achieve near-optimal performance compared to projection for both linear and Lagrangian approaches and for both POD and Greedy compression.

\subsection{Test 2: multi-fidelity training}

\begin{figure}[h!]
\centering
\subfloat[]{
\includegraphics[width=0.45\textwidth]
{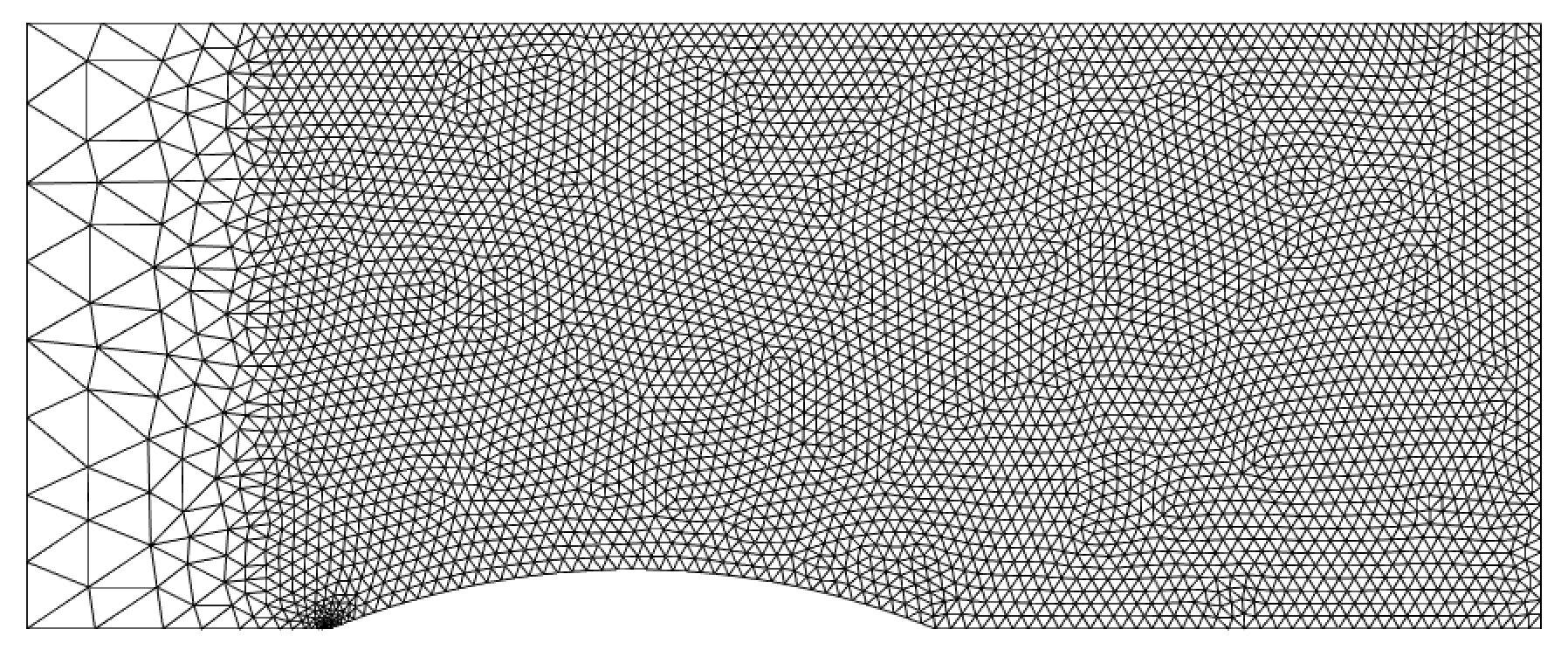}
}
~~~
\subfloat[]{
\includegraphics[width=0.45\textwidth]
{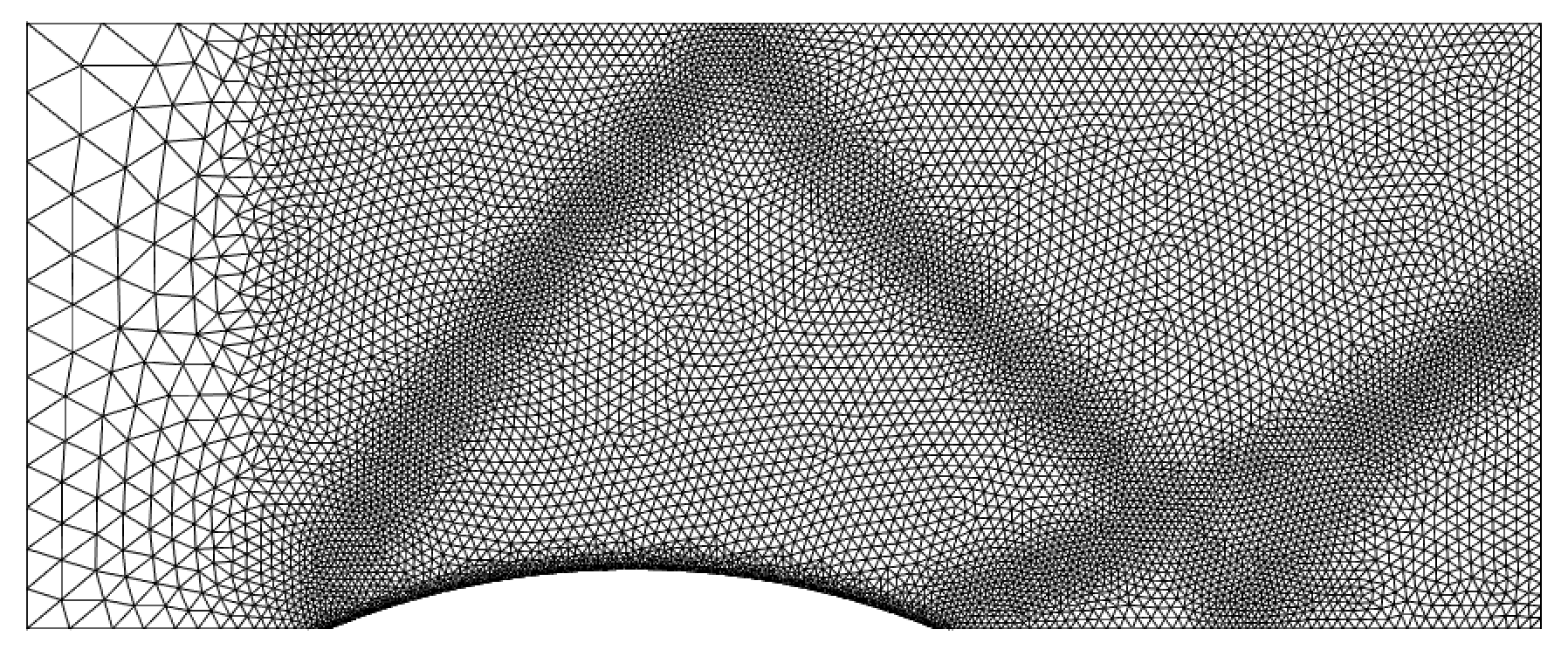}
}
\caption{multi-fidelity training. (a) coarse mesh used for sensor generation. (b) fine mesh used for MOR calculations.
}
\label{fig:test2_meshing}
\end{figure}

We now validate the full offline/online algorithm presented in section \ref{sec:multifidelity}: towards this end, we consider the same hf discretization  and parameter set $\mathcal{P}_{\rm train}$
considered in the previous section to compute the mapping $\Phi$; on the other hand, we use the refined grid depicted in Figure \ref{fig:test2_meshing}(b) with $N_{\rm hf}=402048$ ($N_{\rm e} = 16752$) to generate the hf snapshots. 

As in the previous case,
the mapping $\Phi$ that is obtained applying the registration procedure in Algorithm \ref{alg:registration} consists of three modes ($M=3$);    all three mapping coefficients are well-approximated through  RBF regression. Note that the mapping considered in this test differs from the one in the previous test due to the fact that Algorithm \ref{alg:registration} is fed with a different mesh. Nevertheless, we find that the differences between the two mappings are moderate.

\begin{figure}[h!]
\centering
\subfloat[$\mu_{\rm min}$]{
\includegraphics[width=0.45\textwidth]
{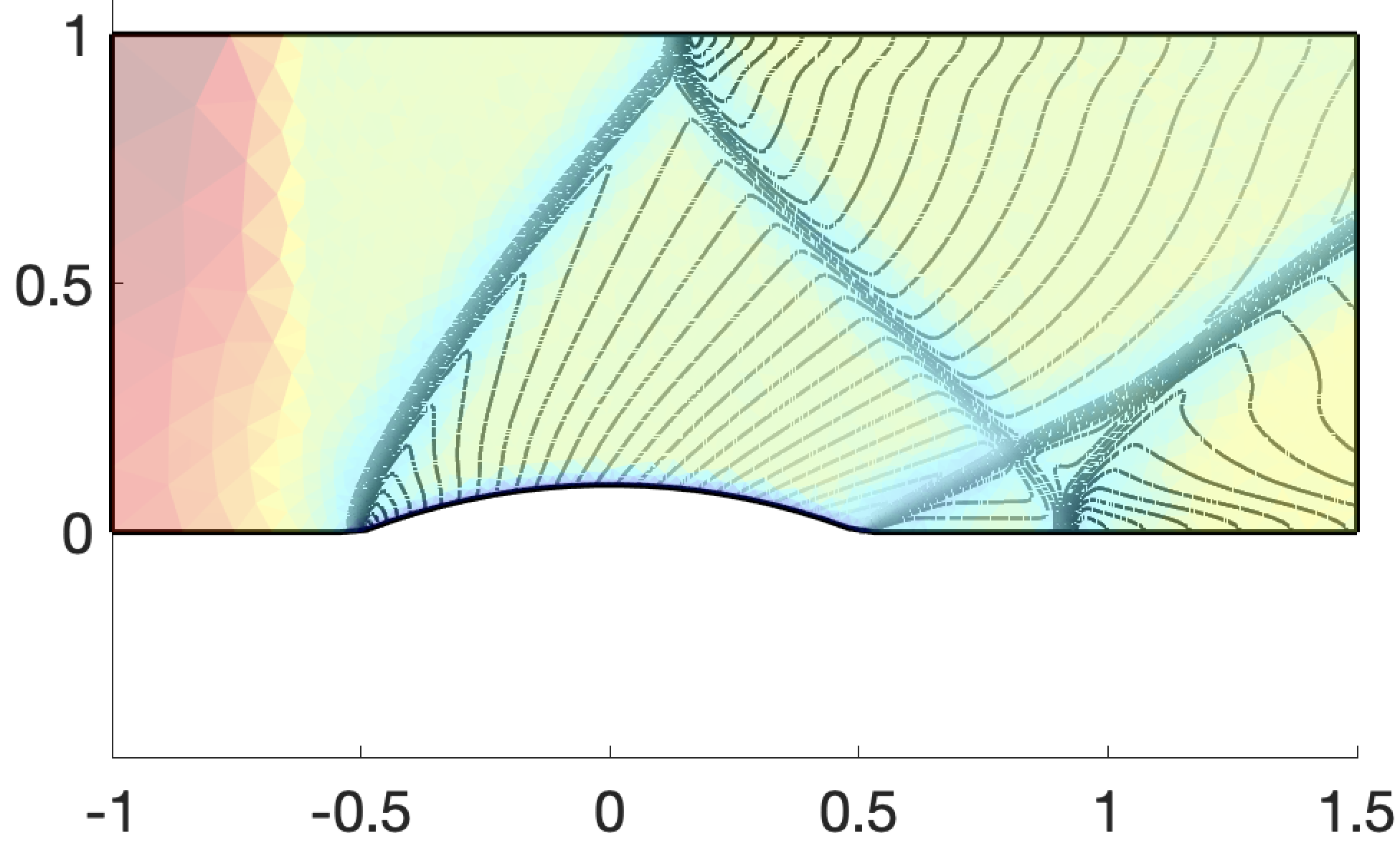}
}
~~~
\subfloat[$\mu_{\rm max}$]{
\includegraphics[width=0.45\textwidth]
{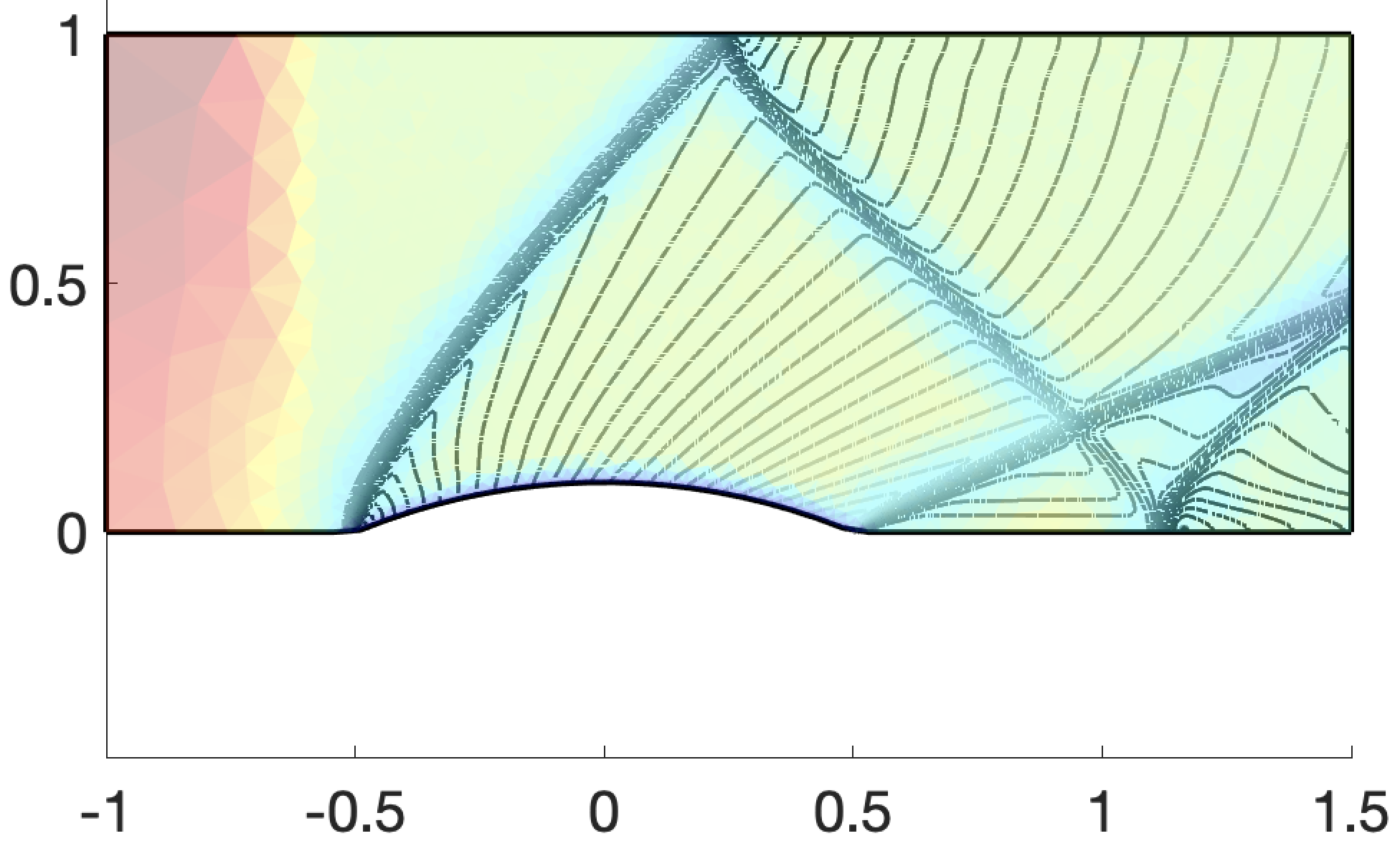}
}

\caption{multi-fidelity training. Comparison of contour lines of Mach number and mesh density $\log_{10}(h)$ for two values of the parameter.
}
\label{fig:test2_meshing_vis}
\end{figure}

In Figure \ref{fig:test2_meshing_vis}, we  investigate the ability of the parametric mesh $\Phi_{\mu}(\mathcal{T}_{\rm hf})$ to track the sharp gradient regions. More in detail, in the background we show 
the mesh density $\log_{10} (h_{\mu})$; in the foreground we show the contour lines of the Mach number, for $\mu_{\rm min} = [0.75,1.7]$ and
$\mu_{\rm max} = [0.8,1.8]$.
Here, the mesh density is defined as 
$h_{\mu}(x) : = \sqrt{ | \texttt{D}_{k,\Phi_{\mu}}   | }$ if $x\in \texttt{D}_{k,\Phi_{\mu}}$, where $\texttt{D}_{k,\Phi_{\mu}} $ is the k-th element of the mesh 
$\Phi_{\mu}(\mathcal{T}_{\rm hf})$.
We observe that the mesh ``follows" the shocks of the solution field: registration is thus able to correctly deform the mesh to track relevant features of the parametric field.

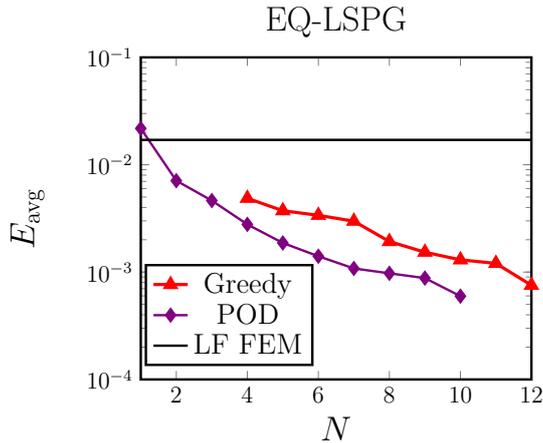
\begin{figure}[h!]
\centering
\begin{tikzpicture}[scale=0.75]
\begin{semilogyaxis}[
title = {\LARGE {EQ-LSPG}},
xlabel = {\LARGE {$N$}},
  ylabel = {\LARGE {$E_{\rm avg}$}},
legend entries = {Greedy, POD, LF FEM},
  line width=1.2pt,
  mark size=3.0pt,
xmin=1,   xmax=12,
 ymin=0.0001,   ymax=0.1,
ylabel style = {font=\large,yshift=5ex},
 yticklabel style = {font=\large,xshift=0ex},
xticklabel style = {font=\large,yshift=0ex},
legend style={at={(0.01,0.2)},anchor=west,font=\Large}
  ]
 
  \addplot[ultra thick,color=red,mark=triangle*]  table {data/euler/test2/greedy_nonlinear.dat};
  
    \addplot[color=violet,mark=diamond*]  table {data/euler/test2/pod_nonlinear.dat};
   
\addplot[color=black,mark=none] table {data/euler/test2/lf_FEM.dat};
  
\end{semilogyaxis}
\end{tikzpicture} 
 
\caption{multi-fidelity training. Performance of EQ-LSPG for POD and weak-Greedy data compression.
}
\label{fig:test2_greedy_training}
\end{figure}

In Figure \ref{fig:test2_greedy_training}, we show performance of EQ-LSPG for 
POD (based on $n_{\rm train}=121$ snapshots) and weak-Greedy data compression; to facilitate interpretation, we further report the average error 
of the coarse solver. We observe that also in this case the ROM is able to provide accurate predictions for extremely moderate values of the ROB size $N$.
In particular, EQ-LSPG with weak-Greedy sampling is able to achieve average out-of-sample errors below $10^{-3}$ with only $N=12$ hf solves.

\section{Conclusions}
\label{sec:conclusions}

In this work, we developed and numerically assessed a multi-fidelity projection- and registration-based MOR procedure for two-dimensional hyperbolic PDEs in presence of shocks. The  key features of our approach are
(i) a general (i.e., independent of the underlying PDE) registration procedure for the computation of the mapping $\Phi$ that tracks moving features of the solution field; 
(ii) an hyper-reduced LSPG ROM for the computation of the solution coefficients; and
(iii) a multi-fidelity approach based on coarse simulations to train the mapping $\Phi$ and   Greedy sampling in parameter,  to reduce  offline costs.
We illustrate the many pieces of our formulation through the vehicle of a supersonic inviscid flow past a bump.

We wish to extend the present work in several directions.
First, our multi-fidelity approach does not include a feedback control on the accuracy of the coarse simulations: for this reason, it might be brittle for more involved problems. It is thus important to devise  robust multi-fidelity strategies that are able to correct the inaccuracies of the coarse simulations.
Second, we wish to relax  the bijectivity-in-$\Omega$ constraint in the registration algorithm by suitably extending the field outside the domain of interest: this would allow to increase the flexibility of   our approach --- particularly, in the presence of fictitious boundaries in the computational domains  --- and ultimately improve performance.
Third, as stated in the introduction, we wish to combine our 
$r$-type, registration-based,  parametric mesh adaptivity technique with $h$-type adaptivity.

\section*{Acknowledgements}
The authors thank  Professor Angelo Iollo (Inria Bordeaux),   Dr. Cédric Goeury and Dr.  Angélique Ponçot (EDF) for fruitful discussions.
The authors acknowledge the support by European Union’s Horizon 2020 research and innovation programme under the Marie Skłodowska-Curie Actions, grant agreement 872442 (ARIA).
Tommaso Taddei also acknowledges the support of IdEx Bordeaux (projet EMERGENCE 2019).

\appendix

\section{Further numerical investigations}
\label{sec:further_numerics}

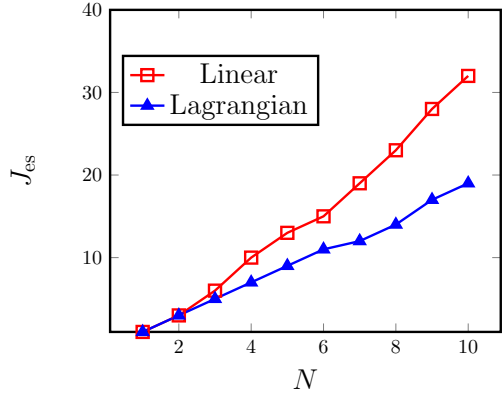
\begin{figure}[h!]
\centering
\begin{tikzpicture}[scale=0.75]
\begin{axis}[
xlabel = {\Large {$N$}},
  ylabel = {\Large {$J_{\rm es}$}},
legend entries = {Linear, Lagrangian},
  line width=1.2pt,
  mark size=3.0pt,
  ylabel style = {font=\large,yshift=2ex},
 ymin=1,   ymax=40,
legend style={at={(0.03,0.75)},anchor=west,font=\Large}
  ]
 
  \addplot[color=red,mark=square]  table {data/euler/test1/lin_test.dat};
  
  \addplot[color=blue,mark=triangle*] table {data/euler/test1/nonlin_test.dat};
    
\end{axis}
\end{tikzpicture}

\caption{single-fidelity training; size of the empirical test space for $tol_{\rm es}=10^{-3}$ for linear and Lagrangian ROMs.
}
\label{fig:test1_empirical_testspace}
\end{figure}

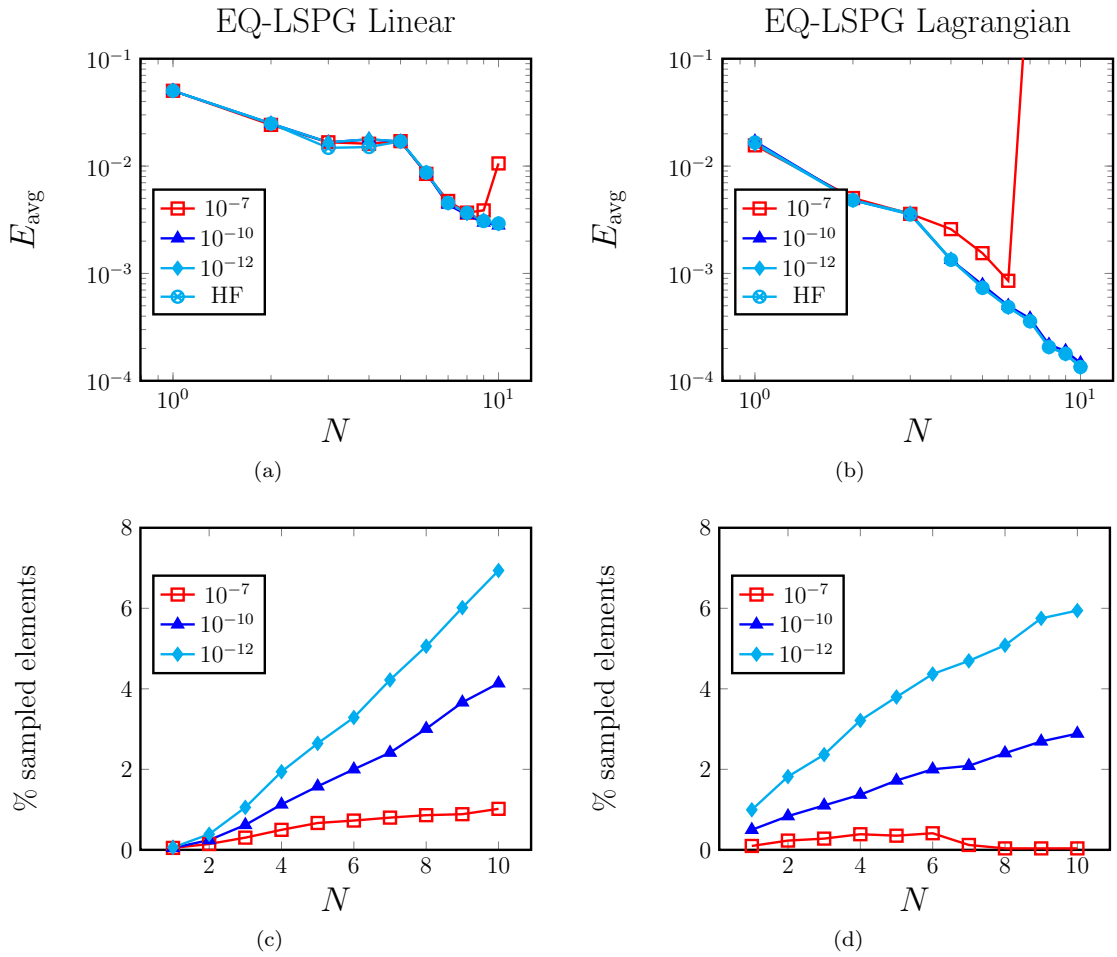
\begin{figure}[h!]
\centering
\subfloat[]{
\begin{tikzpicture}[scale=0.75]
\begin{loglogaxis}[
title = {\LARGE {EQ-LSPG Linear}},
xlabel = {\LARGE {$N$}},
  ylabel = {\LARGE {$E_{\rm avg}$}},
legend entries = {$10^{-7}$,$10^{-10}$,$10^{-12}$,HF},
  line width=1.2pt,
  mark size=3.0pt,
 ymin=0.0001,   ymax=0.1,
ylabel style = {font=\large,yshift=5ex},
 yticklabel style = {font=\large,xshift=0ex},
xticklabel style = {font=\large,yshift=0ex},
legend style={at={(0.03,0.4)},anchor=west,font=\large}
]
\addplot[color=red,mark=square]  table {data/euler/test1/EQ/linear_1m7.dat};
  
  \addplot[color=blue,mark=triangle*] table {data/euler/test1/EQ/linear_1m10.dat};
  
    \addplot[color=cyan,mark=diamond*] table {data/euler/test1/EQ/linear_1m12.dat};
   
    \addplot[color=cyan,mark=otimes] table {data/euler/test1/EQ/linear_hf.dat};   
  
\end{loglogaxis}
\end{tikzpicture}
}
~~~
\subfloat[]{
\begin{tikzpicture}[scale=0.75]
\begin{loglogaxis}[
title = {\LARGE {EQ-LSPG Lagrangian}},
xlabel = {\LARGE {$N$}},
  ylabel = {\LARGE {$E_{\rm avg}$}},
legend entries = {$10^{-7}$,$10^{-10}$,$10^{-12}$,HF},
  line width=1.2pt,
  mark size=3.0pt,
 ymin=0.0001,   ymax=0.1,
ylabel style = {font=\large,yshift=5ex},
 yticklabel style = {font=\large,xshift=0ex},
xticklabel style = {font=\large,yshift=0ex},
legend style={at={(0.03,0.4)},anchor=west,font=\large}
]
\addplot[color=red,mark=square]  table {data/euler/test1/EQ/nonlinear_1m7.dat};
  
  \addplot[color=blue,mark=triangle*] table {data/euler/test1/EQ/nonlinear_1m10.dat};
  
    \addplot[color=cyan,mark=diamond*] table {data/euler/test1/EQ/nonlinear_1m12.dat};
   
    \addplot[color=cyan,mark=otimes] table {data/euler/test1/EQ/nonlinear_hf.dat};   
  
\end{loglogaxis}
\end{tikzpicture}
}

\subfloat[]{
\begin{tikzpicture}[scale=0.75]
\begin{axis}[
xlabel = {\LARGE  {$N$}},
  ylabel = {\Large {$\%$ sampled elements} },
ylabel style = {font=\large,yshift=5ex},
 yticklabel style = {font=\large,xshift=0ex},
xticklabel style = {font=\large,yshift=0ex},
  line width=1.2pt,
  mark size=3.0pt,
  ymin=0,   ymax=8,
legend entries = {$10^{-7}$,$10^{-10}$,$10^{-12}$},
legend style={at={(0.03,0.7)},anchor=west,font=\large}
  ]
 
 \addplot[color=red,mark=square]  table {data/euler/test1/EQ/linear_Qm7.dat};

  \addplot[color=blue,mark=triangle*] table {data/euler/test1/EQ/linear_Qm10.dat};
  
    \addplot[color=cyan,mark=diamond*] table {data/euler/test1/EQ/linear_Qm12.dat};
 
\end{axis}
\end{tikzpicture}
}
~~~
\subfloat[]{
\begin{tikzpicture}[scale=0.75]
\begin{axis}[
xlabel = {\LARGE  {$N$}},
  ylabel = {\Large {$\%$ sampled elements} },
ylabel style = {font=\large,yshift=5ex},
 yticklabel style = {font=\large,xshift=0ex},
xticklabel style = {font=\large,yshift=0ex},
  line width=1.2pt,
  mark size=3.0pt,
  ymin=0,   ymax=8,
legend entries = {$10^{-7}$,$10^{-10}$,$10^{-12}$},
legend style={at={(0.03,0.7)},anchor=west,font=\large}
  ]
 
 \addplot[color=red,mark=square]  table {data/euler/test1/EQ/nonlinear_Qm7.dat};

  \addplot[color=blue,mark=triangle*] table {data/euler/test1/EQ/nonlinear_Qm10.dat};
  
    \addplot[color=cyan,mark=diamond*] table {data/euler/test1/EQ/nonlinear_Qm12.dat};
 
\end{axis}
\end{tikzpicture}
}

\caption{single-fidelity training; hyper-reduction for linear and Lagrangian ROMs.
(a)-(b) behavior of relative error $E_{\rm avg}$ for various tolerances $tol_{\rm eq}$ (cf. \eqref{eq:lsqnonneg}).
(c)-(d) percentage of sampled elements 
$Q/N_{\rm e} \cdot 100$ for the same tolerances $tol_{\rm eq}$.
}
\label{fig:test1_empirical_quadrature}
\end{figure}

\begin{figure}[h!]
\centering
\subfloat[]{
\begin{tikzpicture}[scale=0.75]
\begin{axis}[
title = {\Large {$tol_{\rm eq}=10^{-10}$}},
xlabel = {\Large {$N$}},
  ylabel = {\Large {$\%$ sampled elements}},
legend entries = {fine discr., coarse discr.},
  line width=1.2pt,
  mark size=3.0pt,
  ylabel style = {font=\large,yshift=2ex},
 ymin=0,   ymax=8,
legend style={at={(0.03,0.75)},anchor=west,font=\Large}
  ]
 
  \addplot[color=red,mark=square]  table {data/euler/test2/EQhf_tolm10.dat};
  
  \addplot[color=blue,mark=triangle*] table {data/euler/test2/EQlf_tolm10.dat};
    
\end{axis}
\end{tikzpicture}
}
~~~ 
  \subfloat[]{
 \begin{tikzpicture}[scale=0.75]
\begin{axis}[
title = {\Large {$tol_{\rm eq}=10^{-12}$}},
xlabel = {\Large {$N$}},
  ylabel = {\Large {$\%$ sampled elements}},
legend entries = {fine discr., coarse discr.},
  line width=1.2pt,
  mark size=3.0pt,
  ylabel style = {font=\large,yshift=2ex},
 ymin=0,   ymax=8,
legend style={at={(0.03,0.75)},anchor=west,font=\Large}
  ]
 
  \addplot[color=red,mark=square]  table {data/euler/test2/EQhf_tolm12.dat};
  
  \addplot[color=blue,mark=triangle*] table {data/euler/test2/EQlf_tolm12.dat};
 \end{axis}
\end{tikzpicture}
}
\caption{effect of discretization on hyper-reduction. 
Percentage of sampled elements 
$Q/N_{\rm e} \cdot 100$ for two tolerances $tol_{\rm eq}$ and for fine and coarse discretizations (cf. Figure \ref{fig:test2_meshing}).
}
\label{fig:test2_empirical_quadrature_more}
\end{figure}
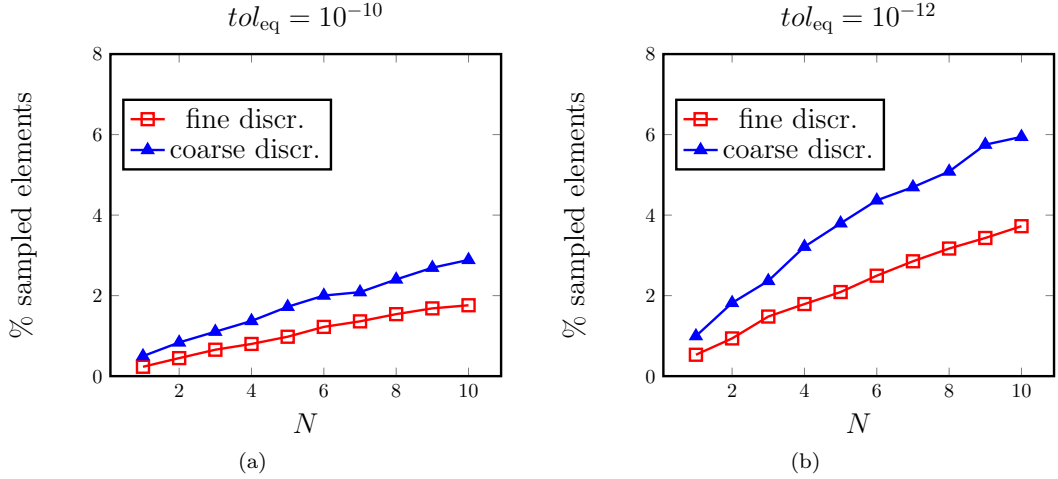

We present here further numerical results to better illustrate the performance of our method. We state upfront that in the results of Figures \ref{fig:test1_empirical_testspace}, \ref{fig:test1_empirical_quadrature}, \ref{fig:test2_empirical_quadrature_more},
we show results for POD data compression.

In  Figure \ref{fig:test1_empirical_testspace} we show the size of the test ROB $\mathbf{Y}$ as obtained using the Algorithm described in section \ref{sec:ROM} for both linear and Lagrangian ROMs. We observe that $J_{\rm es}$ is considerably larger for the linear ROM: registration thus also helps reduce the size of the test space required for stability.

Figure \ref{fig:test1_empirical_quadrature} investigates performance of the hyper-reduction procedure: we show the behavior of the out-of-sample error $E_{\rm avg}$ for different EQ tolerances in \eqref{eq:lsqnonneg}; we further show the percentage of sampled elements $Q/N_{\rm e} \cdot 100$ selected by the EQ procedure. We remark that  EQ ensures accurate performance for $tol_{\rm eq} \leq 10^{-10}$ for all values of $N$ considered and for both linear and Lagrangian ROMs. Interestingly, the linear ROM requires slightly more sampled elements: we conjecture that this is due to the larger size of the test space. 

In Figure \ref{fig:test2_empirical_quadrature_more}, we illustrate the effect of discretization on hyper-reduction: 
we show the percentage of sampled elements $Q/N_{\rm e} \cdot 100$ selected by the EQ procedure for two tolerances, several values of the trial ROB size $N$, and for the two meshes considered in this work (cf. Figure \ref{fig:test2_meshing}). We find that the absolute value of sampled elements weakly depends on the underlying FE  mesh; as a result, hyper-reduction becomes more and more effective as $N_{\rm e}$ increases.

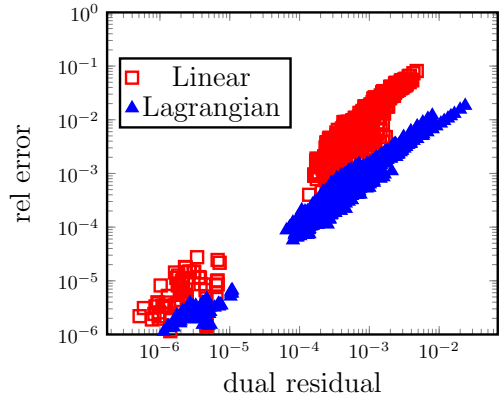
\begin{figure}[h!]
\centering
\begin{tikzpicture}[scale=0.75]
\begin{loglogaxis}[
xlabel = {\Large {dual residual}},
  ylabel = {\Large {rel error}},
legend entries = {Linear, Lagrangian},
  line width=1.2pt,
  mark size=3.0pt,
  ylabel style = {font=\large,yshift=2ex},
 ymin=0.000001,   ymax=1,
legend style={at={(0.03,0.75)},anchor=west,font=\Large}
  ]
 
  \addplot[only marks,color=red,mark=square]  table {data/euler/test1/res_linear.dat};
  
  \addplot[only marks,color=blue,mark=triangle*] table {data/euler/test1/res_nonlinear.dat};
    
\end{loglogaxis}

\end{tikzpicture}

\caption{single-fidelity training; dual residual norm estimation. Comparison between dual residual norm and exact relative error for various ROMs  and  $\mu \in \mathcal{P}_{\rm train}$.
}
\label{fig:test2_dual_residual}
\end{figure}

In Figure \ref{fig:test2_dual_residual}, we investigate the relationship between dual residual \eqref{eq:dual_residual_norm} and relative $L^2$ error for linear and Lagrangian ROMs. More precisely, during each step of the weak-greedy algorithm, we compute both dual residual and relative $L^2$ error for all training points; then, we show the results for all $N=4,\ldots,12$. We observe that there is a strong correlation  between error and dual residual:  this motivates the use of dual residual norm to drive the Greedy algorithm and also as error indicator during the online stage. 
We remark that the points associated with the relative error below $10^{-5}$ correspond to parameters that are sampled by the greedy procedure (see Algorithm \ref{alg:weak_greedy}).

\section{Mesh generation}
\label{sec:mesh_generation}

For completeness, we provide the definition of the mesh size function employed to generate the mesh in Figure \ref{fig:test2_meshing}(b).
We here use the Matlab suite \texttt{distmesh}: we refer to the documentation available at 
\url{persson.berkeley.edu/distmesh/}
for further details.
We envision that the present approach might be greatly improved both in terms of accuracy and in terms of offline computational costs;  the use of state-of-the-art adaptive FE techniques might also be important to automatize the refinement procedure.
Given the coarse simulation $(\mathcal{T}_{\rm hf,c}, \mathbf{U}_{\bar{\mu}}^{\rm hf,c} )$, we define the Mach number ${\rm Ma}^{\rm hf,c} $ and we compute the local averages 
$\mathfrak{s}_1^{\rm c},\ldots,
\mathfrak{s}_{N_{\rm e}^{\rm c}}^{\rm c}$ such that
$$
\mathfrak{s}_k^{\rm c} : =
\int_{ \texttt{D}_{k}^{\rm c}    }
\, \|  \nabla {\rm Ma}^{\rm hf,c}  \|_2^2 \, dx,
\quad
k=1,\ldots,N_{\rm e}^{\rm c}.
$$
We then reorder the elements so that
$\mathfrak{s}_1^{\rm c}\geq  \mathfrak{s}_2^{\rm c}\geq \ldots$; given $n_1=n_2=0.1 \cdot N_{\rm e}^{\rm c}$, we define the 
barycenters $\{  x_j^{\rm c} \}_j$ and the size function 
$$
h^{\rm tmp}(x)
=\min \left\{
3 h_0 + \frac{1}{4} \min\left\{
{\rm dist} 
\left(x, \{  x_j^{\rm c} \}_{j=1}^{n_1} \right), \;
2 h_0 +
{\rm dist} 
\left(x,  \{  x_j^{\rm c} \}_{j=n_1+1}^{n_1+n_2} \right)
\right\},
\;
\bar{h}(x)
\right\}
$$
where $h_0=0.007$, 
$$
\bar{h}(x) = \min \left\{
2 h_0 +{\rm dist}_{\rm bump}(x), \;
6 h_0 + ( -0.6  - x_1  )_+
\right\},
$$
and ${\rm dist}_{\rm bump}(x)$ is the distance of the point $x$ from the semicircular bump. The size function  $h^{\rm tmp}$ measures the proximity to the regions where the gradient of the Mach number is large: it thus leads to mesh refinement in the proximity of the shocks.

The size function  $h^{\rm tmp}$ is excessively irregular for mesh generation purposes: for this reason, we project $h^{\rm tmp}$ over a $100\times 100$ $\texttt{p}=2$ structured uniform grid over $\Omega_{\rm box}=(-1,1.5)\times (0,1)$ and we compute a moving average
with respect to both coordinates; the resulting FE field $h^{\star}$  is passed to the mesh generation routine \texttt{distmesh2d} to generate the  $\texttt{p}=1$ FE grid; 
finally, we perform an iteration of uniform refinement
(see the \texttt{distmesh} routine \texttt{uniref}) to obtain the mesh in Figure \ref{fig:test2_meshing}(b).

\bibliographystyle{abbrv}	
\bibliography{all_refs}

\end{document}